\documentclass[10pt,twoside]{article}
\usepackage{pifont}
\usepackage{bbm}
\usepackage{hyperref}
\usepackage{amsfonts}
\usepackage{amssymb}
\usepackage{amsmath}
\usepackage{amscd}
\usepackage{mathrsfs}
\input amssym.def
\input amssym.tex
\usepackage[all]{xy}
\def\bc{\begin{center}}
\def\ec{\end{center}}
\def\no{\noindent}
\def\a{\alpha}
\def\b{\beta}

\def\c{\cdot}

\def\d{\delta}

\def\D{\Delta}

\def\g{\gamma}

\def\o{\otimes}

\def\v{\varepsilon}

\begin{document}
\begin{center}
{\Large\bf Hopf-Galois algebras and their Poisson structures$^{\ast}$}

\vspace{4mm}

{\bf Huihui Zheng,  Liangyun Zhang}
$$\begin{array}{rllr}
&$~College of Science, Nanjing Agricultural University,
Nanjing 210095, China$
\end{array}$$
\begin{figure}[b]
\rule[-2.5truemm]{5cm}{0.1truemm}\\[2mm]
\no $^{\ast}$This work is supported by Natural Science Foundation (11571173).
\end{figure}
\begin{minipage}{144.6mm}

\vskip1cm \footnotesize{\textbf{Abstract}:
 As is known to all, Hopf-Galois objects have a significant research value for analyzing tensor categories of comodules and classification questions of pointed Hopf algebras, and are natural generalizations of Hopf algebras with a Galois-theoretic flavour. In this paper, we mainly prove a criterion for an Ore extension of a Hopf-Galois algebra to be a Hopf-Galois algebra, and introduce the conception of Poisson Hopf-Galois algebras, and establish the relationship between Poisson Hopf-Galois algebras and Poisson Hopf algebras. Moreover, we study Poisson Hopf-Galois structures on Poisson polynomial algebras, and mainly give a necessary and sufficient condition for the Poisson enveloping algebra of a Poisson Hopf-Galois algebra to be a Hopf-Galois algebra.
\vspace{3mm}}
\vspace{3mm}\\
\textbf{2010 AMS Classification Number}: 16T05, 17B63
\end{minipage}
\end{center}

\vspace{4mm}

\vskip1mm

\begin{center}
{\bf \S1\quad Introduction}
\end{center}

 A (right) Hopf-Galois object \cite{H. F. Kreimer} $R$ over a Hopf algebra $H$ is a (right) comodule algebra $R$ over $H$ with a trivial coinvariants $R^{coH}=k$ (a field) and such that the Galois map
 $$can:R\o R\rightarrow R\o H,~~can(r\o s)=rs_{(0)}\o s_{(1)}, ~\text{for any}~r,s\in R$$
 is bijective, where $\rho(s)=s_{(0)}\o s_{(1)}\in R\o H$ represents the comodule structure of $R$.

 Hopf-Galois objects are important for analyzing tensor categories of comodules \cite{P. Schauenburg 1, K. H. Ulbrich} and classification questions of pointed Hopf algebras \cite{N. Andruskiewitsch, I. Angiono}, and are natural generalizations of Hopf algebras with a Galois-theoretic flavour.

Quantum torsors were introduced by Grunspan \cite{C. Grunspan} with more axioms and simplified later by Schauenburg \cite{P. Schauenburg}. Since there exists a relationship between this quantum torsor and some Hopf-Galois object, the authors call such a structure a Hopf-Galois algebra in \cite{J. Bichon}. The following result relates what we called Hopf-Galois algebras and classical Hopf-Galois objects.

\textbf{Theorem} Let $R$ be an algebra. The following assertion given in \cite{J. Bichon} are equivalent:

(1) There exists an algebra map $\mu: R\rightarrow R\o R^{op}\o R$ making $R$ into a Hopf-Galois algebra.

(2) $R$ is a right Hopf-Galois object over some Hopf algebra.
\vspace{2mm}

 The result ensure that there exists a definition of Hopf-Galois objects that does not use any Hopf algebra, in terms of an algebra map
$\mu: R\rightarrow R\o R^{op}\o R$, subject to certain axioms.

In \cite{J. Bichon}, the authors have proved that eery Hopf algebra is Hopf-Galois algebra, and a Hopf-Galois algebra $R$ arises from a Hopf algebra if and only if there exists an algebra map $\a:R\rightarrow k$. So, Hopf-Galois algebras, like Hopf-Galois objects, are a generalization of Hopf algebras. In addition, Hopf-Galois structures on a generalized ambiskew ring were discussed, and some examples of Hopf-Galois algebras rather than Hopf algebras were provided in \cite{J. Bichon}.

In \cite{J. Lü}, the authors have found that any pointed Poisson Hopf algebra has a Hopf-Galois structure. Naturally, we can consider the Poisson version of a Hopf-Galois algebra. As we expected, in section 3, we establish the correspondence relation between Poisson Hopf algebras and Poisson Hopf-Galois algebras.

Ore extensions of Hopf algebras (HA) and Poisson Hopf algebras (PHA) have been considered in \cite{A. N. Panov,J. Lü}, so, we have a question as follows, that is, how we study Ore extensions of Hopf-Galois algebras (HGA) and Poisson Hopf-Galois algebras (PHGA). In section 2, we give the necessary and sufficient conditions for the Ore extension of a Hopf-Galois algebra to be a Hopf-Galois algebra. Similarly, in section 4, the Poisson Ore extension of a Poisson Hopf-Galois algebra is also studied.\\

\unitlength 1mm 
\linethickness{0.4pt}
\ifx\plotpoint\undefined\newsavebox{\plotpoint}\fi 
\begin{picture}(92.5,26.25)(0,0)
\put(27.25,24.25){HA}
\put(27.25,5){HA}
\put(57.25,24.25){HGA}
\put(57.5,5){HGA}
\put(34.5,25.25){\vector(1,0){21.5}}
\put(34.5,6){\vector(1,0){21.5}}
\put(29.5,22){\vector(0,-1){13.75}}
\put(60.5,22){\vector(0,-1){13.25}}
\put(22.75,15){OE}
\put(75.25,24.25){PHA}
\put(107,24.25){PHGA}
\put(75.25,5){PHA}
\put(107,5){PHGA}
\put(78.75,21.75){\vector(0,-1){13.25}}
\put(85.25,25.5){\vector(1,0){19.75}}
\put(112,21.5){\vector(0,-1){13.25}}
\put(85.25,6){\vector(1,0){19.75}}
\put(73.25,15){OE}
\put(60.5,15){OE $?$}
\put(112.5,15){OE $?$}
\end{picture}

Here OE denotes the Ore extension of an algebra.

In \cite{Sei1}, the author has proved that the Poisson enveloping algebra of a Poisson Hopf algebra has a Hopf algebra structure. Naturally, we can consider whether Poisson enveloping algebra of a Poisson Hopf-Galois algebra has a Hopf-Galois algebra structure. In section 5, we find a necessary and sufficient condition for the Poisson enveloping algebra $U(R)$ of a Poisson Hopf-Galois algebra $R$ to be a Hopf-Galois algebra.\\

\unitlength 1mm 
\linethickness{0.4pt}
\ifx\plotpoint\undefined\newsavebox{\plotpoint}\fi 
\begin{picture}(87.75,22)(0,0)
\put(48.75,21){PHA}
\put(50.75,2.75){HA}
\put(47.75,12.5){UE}
\put(53.5,20){\vector(0,-1){14}}
\put(84.5,21){PHGA}
\put(86.5,2.75){HGA}
\put(90.5,20){\vector(0,-1){14}}
\put(58.25,22){\vector(1,0){24.5}}
\put(58.25,4.5){\vector(1,0){25.5}}
\put(90.75,12.5){UE $?$}
\end{picture}

Here UE denotes the Poisson enveloping algebra of a Poisson algebra.

 \textbf{Notation} Throughout this paper, let $k$ be a fixed field. Unless otherwise specified, linearity, modules and $\otimes$ are all meant over $k$. And we freely use the Hopf algebras terminology introduced in \cite{M. E. Sweedler, H. H. Zheng}. For a coalgebra $C$, we write its comultiplication $\Delta(c)$ with $c_{1}\otimes c_{2}$, for any $c\in C$; for a Hopf-Galois algebra $R$, we denote its structure by $\mu(x)=x_{(1)}\otimes x_{(2)}\o x_{(3)}$, for any $x\in R$, in which we omit the summation symbols for convenience.

\vspace{2mm}

\begin{center}
{\bf \S2\quad Hopf-Galois algebras and their Ore extensions}
\end{center}

In this section, we recall the conception of Hopf-Galois algebras, and mainly prove a criterion for an Ore extension of a Hopf-Galois algebra to be a Hopf-Galois algebra.

\vspace{2mm}

\textbf{Definition 2.1} A Hopf-Galois algebra \cite{J. Bichon} is a non-zero algebra $R$ together with an algebra map
$$
\mu:R\rightarrow R\o R^{op}\o R,
$$
where $R\o R^{op}\o R$ is a tensor product algebra, and such that
\begin{eqnarray*}
&&(\mu\o id\o id)\circ \mu=(id\o id\o \mu)\circ \mu,\\
&&(m\o id)\circ \mu=\eta\o id,~~(id\o m)\circ \mu=id\o \eta,
\end{eqnarray*}
where $m:R\o R\rightarrow R$ and $\eta:k\rightarrow R$ denote the respective multiplication and unit of $R$.
We shall write, for any $r\in R,$
$$
\mu(r)=r_{(1)}\o r_{(2)}\o r_{(3)}.
$$
Hence, for any $r\in R$,
$$
r\o 1=r_{(1)}\o r_{(2)}r_{(3)},~~1\o r=r_{(1)}r_{(2)}\o r_{(3)};\eqno(2.1)
$$
and we may set
$$
r_{(1)}\o r_{(2)}\o r_{(3)}\o r_{(4)}\o r_{(5)}:=\mu(r_{(1)})\o r_{(2)}\o r_{(3)}=r_{(1)}\o r_{(2)}\o \mu(r_{(3)}).\eqno(2.2)
$$

For a Hopf-Galois algebra $R$, we call it a Hopf-Galois commutative algebra if it is commutative. In what follows, we call the above algebra map $\mu$ a Hopf-Galois map.

\vspace{2mm}

\textbf{Remark 2.2} (1) Let $R\stackrel{f}\longrightarrow B\rightarrow 1$ be a short exact sequence of algebras. Assume that $R$ is a Hopf-Galois algebra with the Hopf-Galois map $\mu_R$. Then, $B$ is also a Hopf-Galois algebra with the Hopf-Galois map:
$$
 \mu_B: B\longrightarrow B\otimes B^{op}\otimes B,\ \ b\equiv f(r)\mapsto f(r_{(1)})\otimes f(r_{(2)})\otimes f(r_{(3)}),
 $$
where $\mu_R(b)=r_{(1)}\otimes r_{(2)}\otimes r_{(3)}\in R\otimes R^{op}\otimes R.$

In particular, for any (two-sided) ideal $I$ of the algebra $R$, we know that the quotient algebra $R/I$ is also a Hopf-Galois algebra if $R$ is a Hopf-Galois algebra.

(2) Let $R$ be a Hopf-Galois commutative algebra with a Hopf-Galois map $\mu$. Define
$$\mu'(r)=r_{(3)}\o r_{(2)}\o r_{(1)}
$$
for $r\in R$. Then, $\mu'$ is also a Hopf-Galois map on $R$.

{\bf Proof.} (1) It is straightforward by Definition 2.1.

(2) Denote $\mu'(r)$ by $r_{(1)'}\o r_{(2)'}\o r_{(3)'}$ for $r\in R$. Then, we easily know that
$$\begin{array}{rllr}
r_{(1)'}\o r_{(2)'}r_{(3)'}&=&r_{(3)}\o r_{(2)}r_{(1)}=r\o 1, \\
r_{(1)'}r_{(2)'}\o r_{(3)'}&=&r_{(3)}r_{(2)}\o r_{(1)}=1\o r
\end{array}$$
for any $r\in R$. Moreover
\begin{eqnarray*}
\mu'(r_{(1)'})\o r_{(2)'}\o r_{(3)'}&=&\mu'(r_{(3)})\o r_{(2)}\o r_{(1)}\\
&=&r_{(3)(1)'}\o r_{(3)(2)'}\o r_{(3)(3)'}\o r_{(2)}\o r_{(1)}\\
&=&r_{(3)(3)}\o r_{(3)(2)}\o r_{(3)(1)}\o r_{(2)}\o r_{(1)}\\
&=&r_{(5)}\o r_{(4)}\o r_{(3)}\o r_{(2)}\o r_{(1)},\\
r_{(1)'}\o r_{(2)'}\o \mu'(r_{(3)'})&=&r_{(3)}\o r_{(2)}\o \mu'(r_{(1)})\\
&=&r_{(3)}\o r_{(2)}\o r_{(1)(1)'}\o r_{(1)(2)'}\o r_{(1)(3)'}\\
&=&r_{(3)}\o r_{(2)}\o r_{(1)(3)}\o r_{(1)(2)}\o r_{(1)(1)}\\
&=&r_{(5)}\o r_{(4)}\o r_{(3)}\o r_{(2)}\o r_{(1)},
\end{eqnarray*}
so $\mu'(r_{(1)'})\o r_{(2)'}\o r_{(3)'}=r_{(1)'}\o r_{(2)'}\o \mu'(r_{(3)'})$.

In addition, for any $x,y\in R$, we have
$$\begin{array}{rllr}
\mu'(xy)&=&(xy)_{(1)'}\o (xy)_{(2)'}\o (xy)_{(3)'}=(xy)_{(3)}\o (xy)_{(2)}\o (xy)_{(1)}\\
&=&x_{(3)}y_{(3)}\o x_{(2)}y_{(2)}\o x_{(1)}y_{(1)}\\
&=&\mu'(x)\mu'(y),
\end{array}$$
so $\mu'$ is an algebra map.

\vspace{2mm}

We recall the concept of group-like elements in Hopf-Galois algebras.

\vspace{2mm}

\textbf{Definition 2.3} A group-like element \cite{J. Bichon} in a Hopf-Galois algebra $R$ is an invertible element $g\in R$ such that
$$
\mu(g)=g\o g^{-1}\o g.
$$

In the following, we denote by $G_{0}(R)$ the set of group-like elements in $R$. If $g\in G_{0}(R)$, we write $g\c r:=grg^{-1}$.

\vspace{2mm}

We recall the following result and its proof in \cite{J. Bichon}.

\vspace{2mm}

\textbf{Proposition 2.4} Let $H$ be a Hopf-Galois algebra. Then, there is an algebra map $\a:H\rightarrow k$ if and only if $H$ is a Hopf algebra.

{\bf Proof.} If there is an algebra map $\a:H\rightarrow k$, we can prove that $H$ is a Hopf algebra, whose counit is $\a$ and the other structure maps defined by
$$
\D(x)=\a(x_{(2)})x_{(1)}\o x_{(3)},~~S(x)=\a(x_{(1)}x_{(3)})x_{(2)}.\eqno(2.3)
$$

Conversely, if $H$ is a Hopf algebra, define, for any $x\in R$,
$$
\mu(x)=x_{1}\o S(x_{2})\o x_{3},\eqno(2.4)
$$
then $H$ is a Hopf-Galois algebra.
\hfill $\square$

\vspace*{2mm}

\textbf{Example 2.5}  Assume $H_4=k\{x,g\mid x^2=0,g^2=1,xg=-gx\}$ is the 4-dimensional Sweedler Hopf algebra in \cite{S. Montgomery}. Then, it is easy to check that $H_4$ is a Hopf-Galois algebra with a Hopf-Galois map as folows
\begin{eqnarray*}
&&\mu(x)=x\o 1\o 1-g\o gx\o 1+g\o g\o x;\\
&&\mu(g)=g\o g\o g.
\end{eqnarray*}

Here char$k\neq2$.

\vspace*{2mm}

In what follows, we discuss Hopf-Galois algebra structures on given Ore extension.

Let $A$ be an algebra and $\tau$ an algebra endomorphism of $A$. Let $\delta$ be a $\tau$-derivation of $A$. That is, $\delta$ is a $k$-linear endomorphism $A$ such that $\d(ab)=\tau(a)\d(b)+\d(a)b$ for all $a,b\in A$. The Ore extension $A[z;\tau,\d]$ of the algebra $A$ is an algebra generated by the variable $z$ and the algebra $A$ with the relation
$$
za=\tau(a)z+\d(a)
$$
for all $a\in A$.

If $\{a_i\mid i\in I\}$ is a k-basis of $A$, then, by \cite{J. C. McConnell}, $\{a_iz^j\mid i\in I, j\in \mathbb{N}\}$ is a $k$-basis of the Ore extension $A[z;\tau,\d]$.

Furthermore, assume that $A$ has a Hopf-Galois algebra structure. Then, we can define a Hopf-Galois algebra structure on the Ore extension $A[z;\tau,\d]$.

\vspace*{2mm}

\textbf{Definition 2.6} Let $A$ be a Hopf-Galois algebra and $H=A[z;\tau,\delta]$ an Ore extension of $A$. If $H=A[z;\tau,\delta]$ is a Hopf-Galois algebra with $A$ as a Hopf-Galois subalgebra, then we say $H=A[z;\tau,\delta]$ is a Hopf-Galois Ore extension of $A$.

\vspace{2mm}

\textbf{Proposition 2.7} Let $A$ be a Hopf-Galois algebra and $H=A[z;\tau,\delta]$ an Ore extension of $A$. If $H=A[z;\tau,\delta]$ is a Hopf-Galois Ore extension of $A$, and determined by
$$
\mu(z)=z\o 1\o 1+ g\o g^{-1}\o z-g\o g^{-1}z\o 1,
$$
for some invertible element $g\in A$.
Then $g$ is a group like element of Hopf-Galois algebra $A$.

{\bf Proof.} Assume that there is an extension $\mu:A\rightarrow A\o A^{op} \o A$ such that
$$\mu(z)=z\o 1\o 1+ g\o g^{-1}\o z-g\o g^{-1}z\o 1.
$$

It is obvious that $H$ is a free left $A$-module under the left multiplication with the basis $\{z^i\mid i\geq 0\}$.
By comparing

$(\mu\o id\o id)\mu(z)=z\o 1\o 1\o 1\o 1+g\o g^{-1}\o z\o 1\o 1-g\o g^{-1}z\o 1\o 1\o 1+\mu(g)\o g^{-1}\o z-\mu(g)\o g^{-1}z\o 1$
\\
with

$(id\o id\o \mu)\mu(z)=z\o 1\o 1\o 1\o 1+g\o g^{-1}\o z\o 1\o 1+g\o g^{-1}\o g\o g^{-1}\o z-g\o g^{-1}\o g\o g^{-1}z\o 1-g\o g^{-1}z\o 1\o 1\o 1$,
\\
one gets
$\mu(g)=g\o g^{-1}\o g$.
\hfill $\square$

\vspace*{2mm}

The next theorem gives a sufficient and necessary condition for an Ore extension of a Hopf-Galois algebra to have a Hopf-Galois algebra structure.

\vspace*{2mm}

\textbf{Theorem 2.8} Let $A$ be a Hopf-Galois algebra, $H=A[z;\tau,\delta]$ and $g$ a group like element of Hopf-Galois algebra $A$. If $\tau$ is an algebra automorphism, then $H=A[z;\tau,\delta]$ is a Hopf-Galois Ore extension of $A$, and determined by
$$
\mu(z)=z\o 1\o 1+ g\o g^{-1}\o z-g\o g^{-1}z\o 1,
$$
if and only if the following conditions are satisfied: for all $a\in A$,

(1) $\mu(\tau(a))=\tau(a_{(1)})\o a_{(2)}\o a_{(3)}=\tau(a_{(1)})\o \tau(a_{(2)})\o \tau(a_{(3)})$,

(2) $g\c a_{(1)}\o g\c a_{(2)}\o a_{(3)}=\tau(a_{(1)})\o \tau(a_{(2)})\o a_{(3)}$,

(3) $\d(a_{(1)})\o a_{(2)}\o a_{(3)}+ga_{(1)}\o a_{(2)}g^{-1}\o \d(a_{(3)})+ga_{(1)}\o g^{-1}\d(\tau^{-1}(g\c a_{(2)}))\o a_{(3)}=\mu(\d(a)).$

{\bf Proof.} Assume that $\mu:A\rightarrow A\o A^{op} \o A$ can be extended to $H$ by
$$\mu(z)=z\o 1\o 1+ g\o g^{-1}\o z-g\o g^{-1}z\o 1.$$

Since $za=\tau(a)z+\d(a)$, we have
\begin{eqnarray*}
\mu(za)&=&\mu(z)\mu(a)\\
&=&(z\o 1\o 1+ g\o g^{-1}\o z-g\o g^{-1}z\o 1)(a_{(1)}\o a_{(2)}\o a_{(3)})\\
&=&za_{(1)}\o a_{(2)}\o a_{(3)}+ga_{(1)}\o a_{(2)}g^{-1}\o za_{(3)}-ga_{(1)}\o a_{(2)}g^{-1}z\o a_{(3)}\\
&=&\tau(a_{(1)})z\o a_{(2)}\o a_{(3)}+\d(a_{(1)})\o a_{(2)}\o a_{(3)}+ga_{(1)}\o a_{(2)}g^{-1}\o \tau(a_{(3)})z\\
&+&ga_{(1)}\o a_{(2)}g^{-1}\o \d(a_{(3)})-ga_{(1)}\o a_{(2)}g^{-1}z\o a_{(3)}\\
&=&\tau(a_{(1)})z\o a_{(2)}\o a_{(3)}+\d(a_{(1)})\o a_{(2)}\o a_{(3)}+ga_{(1)}\o a_{(2)}g^{-1}\o \tau(a_{(3)})z\\
&+&ga_{(1)}\o a_{(2)}g^{-1}\o \d(a_{(3)})-ga_{(1)}\o g^{-1}(g\c a_{(2)})z\o a_{(3)}\\
&=&(\tau(a_{(1)})\o a_{(2)}\o a_{(3)})(z\o 1\o 1)+\d(a_{(1)})\o a_{(2)}\o a_{(3)}\\
&+&(g\c a_{(1)}\o g\c a_{(2)}\o \tau(a_{(3)}))(g\o g^{-1}\o z)+ga_{(1)}\o a_{(2)}g^{-1}\o \d(a_{(3)})\\
&-&ga_{(1)}\o g^{-1}z\tau^{-1}(g\c a_{(2)})\o a_{(3)}+ga_{(1)}\o g^{-1}\d(\tau^{-1}(g\c a_{(2)}))\o a_{(3)}\\
&=&(\tau(a_{(1)})\o a_{(2)}\o a_{(3)})(z\o 1\o 1)+\d(a_{(1)})\o a_{(2)}\o a_{(3)}\\
&+&(g\c a_{(1)}\o g\c a_{(2)}\o \tau(a_{(3)}))(g\o g^{-1}\o z)+ga_{(1)}\o a_{(2)}g^{-1}\o \d(a_{(3)})\\
&-&(g\c a_{(1)}\o \tau^{-1}(g\c a_{(2)})\o a_{(3)})(g\o g^{-1}z\o 1)+ga_{(1)}\o g^{-1}\d(\tau^{-1}(g\c a_{(2)}))\o a_{(3)},\\
\mu(za)&=&\mu(\tau(a)z+\d(a))=\mu(\tau(a))\mu(z)+\mu(\d(a))\\
&=&\mu(\tau(a))(z\o 1\o 1)+ \mu(\tau(a))(g\o g^{-1}\o z)-\mu(\tau(a))(g\o g^{-1}z\o 1)+\mu(\d(a)).
\end{eqnarray*}

If $H=A[z;\tau,\delta]$ has a Hopf-Galois algebra structure, so we get
\begin{eqnarray*}
\tau(a)_{(1)}\o \tau(a)_{(2)}\o \tau(a)_{(3)}&=&\tau(a_{(1)})\o a_{(2)}\o a_{(3)}=g\c a_{(1)}\o g\c a_{(2)}\o \tau(a_{(3)}),\\
\tau(a)_{(1)}\o \tau(a)_{(2)}\o \tau(a)_{(3)}&=&g\c a_{(1)}\o \tau^{-1}(g\c a_{(2)})\o a_{(3)},\\
\mu(\d(a))&=&\d(a_{(1)})\o a_{(2)}\o a_{(3)}+ga_{(1)}\o a_{(2)}g^{-1}\o \d(a_{(3)})+ga_{(1)}\\
& & \o g^{-1}\d(\tau^{-1}(g\c a_{(2)}))\o a_{(3)}.
\end{eqnarray*}

The last equation coincides with the condition (3). Combining the first above equation with the second above equation, we can show that
$g\c a_{(1)}\o g\c a_{(2)}\o a_{(3)}=\tau(a_{(1)})\o \tau(a_{(2)})\o a_{(3)}$.

Furthermore, by the above equations, we get

$\mu(\tau(a))=\tau(a_{(1)})\o a_{(2)}\o a_{(3)}=g\c a_{(1)}\o g\c a_{(2)}\o \tau(a_{(3)})=\tau(a_{(1)})\o \tau(a_{(2)})\o \tau(a_{(3)}).$

Conversely, if the conditions (1)$\sim$(3) hold, then, $H$ has a Hopf-Galois structure by the previous computation.

\vspace*{2mm}

\begin{center}
{\bf \S3\quad Poisson Hopf-Galois algebras}
\end{center}

In this section, we introduce the conception of Poisson Hopf-Galois algebra, and establish the relationship between Poisson Hopf-Galois algebras and Poisson Hopf algebras.

\vspace*{2mm}

\textbf{Definition 3.1} A Poisson algebra in \cite{L. I. Korogodski} is a commutative algebra $A$ equipped with a bilinear map $\{\cdot,\cdot\}:A\o A\rightarrow A$ such that $A$ is a Lie algebra under $\{\cdot,\cdot\}$ and
$$
\{a,bc\}=\{a,b\}c+b\{a,c\}\eqno(3.1)
$$
for any $a,b,c\in A$.

\vspace{2mm}

\textbf{Remark 3.2} If $A$ and $B$ are Poisson algebras, then, the tensor product algebra $A\o B$ has a Poisson structure defined by
$$
\{a\o b,a'\o b'\}_{A\otimes B}=aa'\o \{b,b'\}+\{a,a'\}\o bb'\eqno(3.2)
$$
 for all $a,a'\in A$, and $b,b'\in B$.

\vspace{2mm}

\textbf{Definition 3.3} A Poisson algebra $A$ is said to be a Poisson Hopf algebra in \cite{L. I. Korogodski} if $A$ is also a Hopf algebra $(A,1,m,\D,\v,S)$ such that
$$
\D(\{a,b\})=\{\D(a),\D(b)\}_{A\o A}\eqno(3.3)
$$
for all $a,b\in A$.
\vspace{2mm}

For Poisson algebras $A$ and $B$, an algebra homomorphism $\phi: A\rightarrow B$ is said to be a Poisson homomorphism, anti-homomorphism) if $\phi$ satisfies the rule
$$
\phi(\{a,b\})=\{\phi(a),\phi(b)\}~~(\phi(\{a,b\})=\{\phi(b),\phi(a)\})
$$
for all $a,b\in A$.

\vspace{2mm}

\textbf{Lemma 3.4} Let $A$ be a Poisson Hopf algebra. Then, by \cite{J. Lü}, $\v(\{a,b\})=0$ and the antipode $S$ is a Poisson anti-automorphism.

\vspace{2mm}

In what follows, we introduce the conception of a Poisson Hopf-Galois algebra by the same axioms like Grunspan's Poisson quantum torsors, but without endomorphism $\theta$ in \cite{C. Grunspan}.

\vspace{2mm}

\textbf{Definition 3.5}  A Poisson Hopf-Galois algebra $R$ is a Poisson algebra and a Hopf-Galois algebra, such that the Hopf-Galois map $\mu$ on $R$ is a Poisson homomorphism, that is, $\mu$ is an algebra map, and for any $a,b\in R,$
$$\mu(\{a,b\})=\{\mu(a),\mu(b)\}.\eqno(3.4)
$$

Note that the Poisson structure on $R\o R^{op}\o R$ is given by
$$
\{x\o y\o z, x'\o y'\o z'\}=\{x,x'\}\o yy'\o zz'-xx'\o \{y,y'\}\o zz'+xx'\o yy'\o \{z,z'\},\eqno(3.5)
$$
for all $x,y,z,x',y',z'\in R$. So, the equality (3.4) can be written into
\begin{eqnarray*}
\{a,b\}_{(1)}\otimes\{a,b\}_{(2)}\otimes\{a,b\}_{(3)}&=&\{a_{(1)},b_{(1)}\}\otimes a_{(2)}b_{(2)}\otimes a_{(3)}b_{(3)}\\
&-&a_{(1)}b_{(1)}\otimes\{a_{(2)},b_{(2)}\}\otimes a_{(3)}b_{(3)}+a_{(1)}b_{(1)}\otimes a_{(2)}b_{(2)}\otimes\{a_{(3)},b_{(3)}\}.
 \end{eqnarray*}

\vspace{2mm}

\textbf{Remark 3.6 } (1) Assume that $R\stackrel{f}\longrightarrow B\rightarrow 1$ is a short exact sequence of Poisson algebras. If $R$ is a Poisson Hopf-Galois algebra with the Hopf-Galois map $\mu_R$, then, $B$ is also a Poisson Hopf-Galois algebra with the Hopf-Galois map:
$$
 \mu_B: B\longrightarrow B\otimes B^{op}\otimes B,\ \ b\equiv f(r)\mapsto f(r_{(1)})\otimes f(r_{(2)})\otimes f(r_{(3)}),
 $$
where $\mu_R(b)=r_{(1)}\otimes r_{(2)}\otimes r_{(3)}\in R\otimes R^{op}\otimes R.$

In particular, for any ideal $I$ of the Poisson algebra $R$ (that is, $RI\subseteq I$ and $\{R, I\}\subseteq I$), we know that the quotient algebra $R/I$ is also a Poisson Hopf-Galois algebra if $R$ is a Poisson Hopf-Galois algebra.

(2) Let $R$ be a Poisson Hopf-Galois algebra with a Hopf-Galois map $\mu$. Define
$$\mu'(r)=r_{(3)}\o r_{(2)}\o r_{(1)}
$$
for $r\in R$. Then, $(R,\mu')$ is also a Poisson Hopf-Galois algebra with the Hopf-Galois map $\mu^\prime$.

{\bf Proof.} (1) On the basis of Remark 2.2, we know that $(B, \mu_B)$ is a Hopf-Galois algebra with the Hopf-Galois map $\mu_B$.

So, by Definition 3.5, we only need to prove that $\mu_B(\{b,b'\})=\{\mu_B(b),\mu_B(b')\}$ for any $b,b'\in B$ in order to prove that $(B, \mu_B)$ is a Poisson Hopf-Galois algebra.

In fact, since there are elements $r,r'\in R$ such that $b=f(r)$ and $b'=f(r')$, we get
\begin{eqnarray*}
\mu_B(\{b,b'\})&=&\mu_B(\{f(r),f(r')\})=\mu_B f(\{r,r'\})\\
&=&f(\{r,r'\}_{(1)})\otimes f(\{r,r'\}_{(2)})\otimes f(\{r,r'\}_{(3)})\\
&\stackrel{(3.4)}=&f(\{r_{(1)},r'_{(1)}\})\otimes f(r_{(2)}r'_{(2)})\otimes f(r_{(3)}r'_{(3)})-f(r_{(1)}r'_{(1)})\otimes f(\{r_{(2)},r'_{(2)}\})\otimes f(r_{(3)}r'_{(3)})\\
&+&f(r_{(1)}r'_{(1)})\otimes f(r_{(2)}r'_{(2)})\otimes f(\{r_{(3)},r'_{(3)}\}),\\
\{\mu_B(b),\mu_B(b')\}&=&\{\mu_B(f(r)),\mu_B(f(r'))\}\\
&=&\{f(r_{(1)})\otimes f(r_{(2)})\otimes f(r_{(3)}),f(r'_{(1)})\otimes f(r'_{(2)})\otimes f(r'_{(3)})\}\\
&\stackrel{(3.5)}=&\{f(r_{(1)}),f(r'_{(1)})\}\otimes f(r_{(2)})f(r'_{(2)})\otimes f(r_{(3)})f(r'_{(3)})\\
&-&f(r_{(1)})f(r'_{(1)})\otimes \{f(r_{(2)}),f(r'_{(2)})\}\otimes f(r_{(3)})f(r'_{(3)})\\
&+&f(r_{(1)})f(r'_{(1)})\otimes f(r_{(2)})f(r'_{(2)})\otimes \{f(r_{(3)}),f(r'_{(3)})\}\\
&=&f(\{r_{(1)},r'_{(1)}\})\otimes f(r_{(2)}r'_{(2)})\otimes f(r_{(3)}r'_{(3)})-f(r_{(1)}r'_{(1)})\otimes f(\{r_{(2)},r'_{(2)}\})\otimes f(r_{(3)}r'_{(3)})\\
&+&f(r_{(1)}r'_{(1)})\otimes f(r_{(2)}r'_{(2)})\otimes f(\{r_{(3)},r'_{(3)}\}),
\end{eqnarray*}
so $\mu_B(\{b,b'\})=\{\mu_B(b),\mu_B(b')\}$.

(2) It is straightforward.

\vspace{2mm}

In what follows, we establish the relationship between Poisson Hopf-Galois algebras and Poisson Hopf algebras.

\vspace{2mm}

\textbf{Proposition 3.7} Let $H$ be a Poisson algebra. Then the following conclusions hold.

(1) Suppose that $H$ is a Hopf-Galois algebra with the Hopf-Galois map $\mu$. If there is an algebra map $\a:H\rightarrow k$ such that $\{H,H\}\subseteq$Ker$\a$, then, $H$ has a structure of Poisson Hopf algebra.

(2) If $H$ is a Poisson Hopf algebra, then, $H$ is also a Poisson Hopf-Galois algebra.

{\bf Proof.} (1) By Proposition 2.3, we know that $H$ is a Hopf algebra whose counit is $\a$ and the other structure maps still defined by Eq.(2.3).

So, in order to prove that $H$ has a structure of Poisson Hopf algebra, we need to prove that $\D(\{x,y\})=\{\D(x),\D(y)\}$ for any $x,y\in H$.

In fact, since $H$ is a Poisson Hopf-Galois algebra, we have $\mu(\{x,y\})=\{\mu(x),\mu(y)\}$, that is, Eq.(3.4) holds.

Hence, we get
\begin{eqnarray*}
\D(\{x,y\})&=&\a(\{x,y\}_{(2)})\{x,y\}_{(1)}\o \{x,y\}_{(3)}\\
&=&\a(x_{(2)}y_{(2)})\{x_{(1)},y_{(1)}\}\o x_{(3)}y_{(3)}-\a(\{x_{(2)},y_{(2)}\})x_{(1)}y_{(1)}\o x_{(3)}y_{(3)}\\
&+&\a(x_{(2)}y_{(2)})x_{(1)}y_{(1)}\o \{x_{(3)},y_{(3)}\}\\
&=&\a(x_{(2)}y_{(2)})\{x_{(1)},y_{(1)}\}\o x_{(3)}y_{(3)}+\a(x_{(2)}y_{(2)})x_{(1)}y_{(1)}\o \{x_{(3)},y_{(3)}\},\\
\{\D(x),\D(y)\}&=&x_1y_1\o \{x_2,y_2\}+\{x_1,y_1\}\o x_2y_2\\
&=&\a(x_{(2)})x_{(1)}\a(y_{(2)})y_{(1)}\o \{x_{(3)},y_{(3)}\}+\{\a(x_{(2)})x_{(1)},\a(y_{(2)})y_{(1)}\}\o x_{(3)}y_{(3)}\\
&=&\a(x_{(2)}y_{(2)})x_{(1)}y_{(1)}\o \{x_{(3)},y_{(3)}\}+\a(x_{(2)}y_{(2)})\{x_{(1)},y_{(1)}\}\o x_{(3)}y_{(3)}.
\end{eqnarray*}

So, $\D(\{x,y\})=\{\D(x),\D(y)\}$, and $H$ is a Poisson Hopf algebra.

(2) Assume that $H$ is a Poisson Hopf algebra. Then, there exist an algebra map $\a=\varepsilon: H\rightarrow k$ such that $\{H,H\}\subseteq$Ker$\a$ by Lemma 3.4.

Define the map $\mu$ by Eq.(2.4), we get $H$ is a Hopf-Galois algebra by Proposition 2.3. In what follows, we only need to check $\mu(\{x,y\})=\{\mu(x),\mu(y)\}$ in order to prove that $H$ is a Poisson Hopf-Galois algebra.

As a matter of fact, since $\D(\{x,y\})=\{\D(x),\D(y)\}$, we have
$$\{x,y\}_1\o \{x,y\}_2\o \{x,y\}_3=x_1y_1\o x_2y_2\o \{x_3, y_3\}+\{x_1, y_1\}\o x_2y_2\o x_3y_3+x_1y_1\o \{x_2, y_2\}\o x_3y_3
.$$

Hence we get
\begin{eqnarray*}
\mu(\{x,y\})&=&\{x,y\}_1\o S(\{x,y\}_2)\o \{x,y\}_3\\
&=&x_1y_1\o S(x_2y_2)\o \{x_3, y_3\}+\{x_1, y_1\}\o S(x_2y_2)\o x_3y_3+x_1y_1\o S(\{x_2, y_2\})\o x_3y_3\\
&=&x_1y_1\o S(x_2y_2)\o \{x_3, y_3\}+\{x_1, y_1\}\o S(x_2y_2)\o x_3y_3\\
&-&x_1y_1\o \{S(x_2),S(y _2)\}\o x_3y_3,\\
\{\mu(x),\mu(y)\}&=&\{x_{(1)}\o x_{(2)}\o x_{(3)}, y_{(1)}\o y_{(2)}\o y_{(3)}\}\\
&=&\{x_1\o S(x_2)\o x_3, y_1\o S(y_2)\o y_3\}\\
&=&\{x_1, y_1\}\o S(x_2)S(y_2)\o x_3y_3-x_1y_1\o \{S(x_2),S(y _2)\}\o x_3y_3\\
&+&x_1y_1\o S(x_2)S(y_2)\o \{x_3, y_3\}\\
&=&\{x_1, y_1\}\o S(x_2y_2)\o x_3y_3-x_1y_1\o \{S(x_2),S(y _2)\}\o x_3y_3\\
&+&x_1y_1\o S(x_2y_2)\o \{x_3, y_3\}.
\end{eqnarray*}

So, $\mu(\{x,y\})=\{\mu(x),\mu(y)\}$, and $H$ is a Poisson Hopf-Galois algebra. $\Box$

\vspace{2mm}

\textbf{Example 3.8}  Let $B=k[g,g^{-1},x]$ be the polynomial algebra. Define a bracket product on $B$ as follows
$$
\{x,g\}=\lambda g x,
$$
where $0\neq\lambda\in k$. Then, there is a Poisson Hopf-Galois algebra on $B$ with the Hopf-Galois map
 $$\mu(g)=g\o g^{-1}\o g,
 $$
 $$
 \mu(x)=1\o 1\o x-1\o xg\o g^{-1}+x\o g\o g^{-1}.
 $$

{\bf Proof.} Assume that $(B, \{\cdot,\cdot\})$ is a Poisson algebra. Then, $\{g,g^{-1}\}=0$ by the equality
$$0=\{g,g^{-1}g\}=\{g,g^{-1}\}g+g^{-1}\{g,g\}.
$$

Since
$$\{x,g^{-1}g\}=\{x,g^{-1}\}g+g^{-1}\{x,g\},
$$
we get
$$\{x,g^{-1}\}=-\lambda g^{-1}x.$$

In what follows, we check that $B$ is a Hopf-Galois algebra.
Because $\mu(g)=g\o g^{-1}\o g$, we get
$$\mu(g^{-1})=g^{-1}\o g\o g^{-1}.$$

By the above discussion, we easily get that
$$\begin{array}{rllr}
g\o 1&=g_{(1)}\o g_{(2)}g_{(3)},\ \ 1\o g=g_{(1)}g_{(2)}\o g_{(3)},\\
g^{-1}\o 1&=(g^{-1})_{(1)}\o (g^{-1})_{(2)}(g^{-1})_{(3)},\ \ 1\o g^{-1}=(g^{-1})_{(1)}(g^{-1})_{(2)}\o (g^{-1})_{(3)},\\
x_{(1)}\o x_{(2)}x_{(3)}&=1\o x-1\o xgg^{-1}+x\o 1=x\o 1,\\
x_{(1)}x_{(2)}\o x_{(3)}&=1\o x-xg\o g^{-1}+xg\o g^{-1}=1\o x,\\
\mu(g_{(1)})\o g_{(2)}\o g_{(3)}&=g_{(1)}\o g_{(2)}\o \mu(g_{(3)}),\\
\mu(x_{(1)})\o x_{(2)}\o x_{(3)}&=1\o 1\o 1\o 1\o x-1\o 1\o 1\o xg\o g^{-1}+1\o 1\o x\o g\o g^{-1}\\
&=x_{(1)}\o x_{(2)}\o \mu(x_{(3)}),
\end{array}$$
so, $B$ is a Hopf-Galois algebra with the Hopf-Galois map $\mu$.

In the following, we prove that $\mu$ is a Poisson map in order to explain $B$ is a Poisson Hopf-Galois algebra.

As a matter of fact, we have no difficulty to prove
$$\begin{array}{rllr}
\{\mu(g),\mu(g^{-1})\}&=&\{g\o g^{-1}\o g,g^{-1}\o g\o g^{-1}\}=0\\
&=&\mu(\{g,g^{-1}\}),\\
\{\mu(x),\mu(g)\}&=&\{1\o 1\o x-1\o xg\o g^{-1}+x\o g\o g^{-1},g\o g^{-1}\o g\}\\
&=&\{1\o 1\o x,g\o g^{-1}\o g\}-\{1\o xg\o g^{-1},g\o g^{-1}\o g\}\\
&+&\{x\o g\o g^{-1},g\o g^{-1}\o g\}\\
&=&g\o g^{-1}\o \{x,g\}+g\o \{xg,g^{-1}\}\o 1+g\o x\o \{g^{-1},g\}\\
&&+\{x,g\}\o 1\o 1-xg\o \{g,g^{-1}\}\o 1+xg\o 1\o \{g^{-1},g\}\\
\end{array}$$
$$\begin{array}{rllr}
&=&g\o g^{-1}\o \{x,g\}+g\o \{xg,g^{-1}\}\o 1+\{x,g\}\o 1\o 1\\
&=&g\o g^{-1}\o \{x,g\}-g\o \{g^{-1},x\}g\o 1+\{x,g\}\o 1\o 1\\
&=&g\o g^{-1}\o \lambda gx-g\o \lambda x\o 1+\lambda gx \o 1\o 1\\
&=&\lambda (g\o g^{-1}\o gx-g\o x\o 1+gx\o 1\o 1)\\
&=&\lambda (g\o g^{-1}\o g)(1\o 1\o x-1\o xg\o g^{-1}+x\o g\o g^{-1})\\
&=&\lambda\mu(g)\mu(x)=\mu(\{x,g\}).
\end{array}$$

In a similar way, we can prove $\mu(\{g^{-1},x\})=\{\mu(g^{-1}),\mu(x)\}$.

\vspace{2mm}

\begin{center}
{\bf \S4\quad Poisson Hopf-Galois structures on polynomial algebras}
\end{center}

\vspace{2mm}

The Poisson polynomial algebra introduced in \cite{Sei3} is a Poisson version of the Ore extension. In this section, we study the Poisson Hopf-Galois structure on a Poisson polynomial algebra.

\vspace{2mm}

\textbf{Definition 4.1} Let $B$ be a Poisson algebra. A Poisson derivation on $B$ is a $k$-linear map $\a$ from $B$ to $B$ which is a derivation with respect to both multiplication and the Poisson bracket, that is, for any $a,b\in B$,
$$\a(ab)=\a(a)b+a\a(b),\eqno(4.1) $$
$$\a(\{a,b\})=\{\a(a),b\}+\{a,\a(b)\}. \eqno(4.2)
$$

\textbf{Remark 4.2} Let $B$ be a Poisson algebra and $\a$ a Poisson derivation on $B$. Suppose that $\delta$ is a derivation with respect to multiplication on $B$, such that for any $a,b\in B$,
$$
\delta(\{a,b\})=\{\delta(a), b\}+\{a, \delta(b)\}+\a(a)\delta(b)-\delta(a)\a(b). \eqno(4.3)
$$

Then, according to \cite{Sei3}, we know that the Poisson structure on $B$ can be extended uniquely to a Poisson structure on the polynomial algebra $R=B[x]$ via
$$
\{x,b\}=\a(b)x+\delta(b) \eqno(4.4)
$$
for any $b\in B$.

In the following, we write $R=B[x;\a,\delta]_{p}$ for this Poisson algebra $B$ and call it a Poisson Ore extension of $B$.

\vspace{2mm}

\textbf{Definition 4.3} Let $B$ be a Poisson Hopf-Galois algebra with the Hopf-Galois map $\mu$, and $B[x;\a,\delta]_{p}$ a Poisson Ore extension of $B$. If $B[x;\a,\delta]_{p}$ is a Poisson Hopf-Galois algebra with $B$ as a Poisson Hopf-Galois subalgebra, then we call  $B[x;\a,\delta]_{P}$ a Poisson Hopf-Galois Ore extension of $B$.

\vspace{2mm}

Naturally, we have the following question: when is $B[x;\a,\delta]_{p}$ a Poisson Hopf-Galois Ore extension of $B$ ?

In what follows, we always assume that $B[x;\a,\delta]_{p}$ is a Poisson Hopf-Galois algebra with $B$ as a Poisson Hopf-Galois subalgebra, whose Hopf-Galois map is given by
$$\mu(x)=x\o 1\o 1-g\o g^{-1}x\o 1+g\o g^{-1}\o x,\eqno(4.5)
$$
where $g$ is a group-like element in $B$.

\vspace{2mm}

\textbf{Theorem 4.4} Let $B$ be a Poisson Hopf-Galois algebra with a group-like element $g$. Then $B[x;\a,\delta]_{p}$ is a Poisson Hopf-Galois Ore extension of $B$, if and only if

(1) the equalities (4.6)$\sim$ (4.8) hold for any $a,b\in B$,
$$\a(b)=g^{-1}\{g,b\},\eqno(4.6) $$
$$
\mu(\a(b))=\a(b_{(1)})\o b_{(2)}\o b_{(3)},\eqno(4.7)
$$
$$
b_{(1)}\o b_{(2)}\o \a(b_{(3)})=b_{(1)}\o g\{g^{-1},b_{(2)}\}\o b_{(3)},\eqno(4.8)
$$
$$
\{g^{-1}, b\}\a(a)=\{g^{-1}, a\}\a(b),\eqno(4.9)
$$

(2) the given map $\delta$ satisfies
$$
\mu(\delta(b))=\delta(b_{(1)})\o b_{(2)}\o b_{(3)}+gb_{(1)}\o \{g^{-1}x,b_{(2)}\}\o b_{(3)}+gb_{(1)}\o g^{-1}b_{(2)}\o \delta(b_{(3)})
\eqno(4.10)
$$
for any $b\in B$.

{\bf Proof.} For any $b\in B$, we have
\begin{eqnarray*}
\mu(\{x,b\})&=&\mu(\a(b)x+\delta(b))=\mu(\a(b))\mu(x)+\mu(\delta(b))\\
&\stackrel{(4.5)}=&\mu(\a(b))(x\o 1\o 1)-\mu(\a(b))(g\o g^{-1}x\o 1)+\mu(\a(b))(g\o g^{-1}\o x)+\mu(\delta(b)),\\
\{\mu(x),\mu(b)\}&=&\{x\o 1\o 1-g\o g^{-1}x\o 1+g\o g^{-1}\o x, b_{(1)}\o b_{(2)}\o b_{(3)}\}\\
&=&\{x\o 1\o 1, b_{(1)}\o b_{(2)}\o b_{(3)}\}-\{g\o g^{-1}x\o 1, b_{(1)}\o b_{(2)}\o b_{(3)}\}\\
&+&\{g\o g^{-1}\o x, b_{(1)}\o b_{(2)}\o b_{(3)}\}\\
&\stackrel{(3.5)}=&\{x, b_{(1)}\}\o b_{(2)}\o b_{(3)}-\{g, b_{(1)}\}\o g^{-1}xb_{(2)}\o b_{(3)}+gb_{(1)}\o \{g^{-1}x,b_{(2)}\}\o b_{(3)}\\
&+&\{g, b_{(1)}\}\o g^{-1}b_{(2)}\o xb_{(3)}-gb_{(1)}\o \{g^{-1}, b_{(2)}\}\o xb_{(3)}\\
&+&gb_{(1)}\o g^{-1}b_{(2)}\o \{x, b_{(3)}\}\\
&\stackrel{(4.4)}=&(\a(b_{(1)})x+\delta(b_{(1)}))\o b_{(2)}\o b_{(3)}-\{g, b_{(1)}\}\o g^{-1}xb_{(2)}\o b_{(3)}\\
&+&gb_{(1)}\o \{g^{-1}x,b_{(2)}\}\o b_{(3)}+\{g, b_{(1)}\}\o g^{-1}b_{(2)}\o xb_{(3)}\\
&-&gb_{(1)}\o \{g^{-1}, b_{(2)}\}\o xb_{(3)}+gb_{(1)}\o g^{-1}b_{(2)}\o (\a(b_{(3)})x+\delta(b_{(3)}))\\
&=&\a(b_{(1)})x\o b_{(2)}\o b_{(3)}+\delta(b_{(1)})\o b_{(2)}\o b_{(3)}-\{g, b_{(1)}\}\o g^{-1}xb_{(2)}\o b_{(3)}\\
&+&gb_{(1)}\o \{g^{-1}x,b_{(2)}\}\o b_{(3)}+\{g, b_{(1)}\}\o g^{-1}b_{(2)}\o xb_{(3)}\\
&-&gb_{(1)}\o \{g^{-1}, b_{(2)}\}\o xb_{(3)}+gb_{(1)}\o g^{-1}b_{(2)}\o \a(b_{(3)})x\\
&+&gb_{(1)}\o g^{-1}b_{(2)}\o\delta(b_{(3)})\\
&=&(\a(b_{(1)})\o b_{(2)}\o b_{(3)})(x\o 1\o 1)+\delta(b_{(1)})\o b_{(2)}\o b_{(3)}\\
&-&(g^{-1}\{g, b_{(1)}\}\o b_{(2)}\o b_{(3)})(g\o g^{-1}x\o 1)+gb_{(1)}\o \{g^{-1}x,b_{(2)}\}\o b_{(3)}\\
&+&(g^{-1}\{g, b_{(1)}\}\o b_{(2)}\o b_{(3)})(g\o g^{-1}\o x)\\
&-&(b_{(1)}\o g\{g^{-1}, b_{(2)}\}\o b_{(3)})(g\o g^{-1}\o x)\\
&+&(b_{(1)}\o b_{(2)}\o \a(b_{(3)}))(g\o g^{-1}\o x)+gb_{(1)}\o g^{-1}b_{(2)}\o\delta(b_{(3)}).
\end{eqnarray*}

Assume that $B[x;\a,\delta]_{p}$ is a Poisson Hopf-Galois algebra. Then, $\mu(\{x,b\})=\{\mu(x),\mu(b)\})$. Thus, by the above computation, we get
\begin{eqnarray*}
&&\mu(\a(b))=\a(b_{(1)})\o b_{(2)}\o b_{(3)},\\
&&\mu(\a(b))=g^{-1}\{g, b_{(1)}\}\o b_{(2)}\o b_{(3)},\\
&&\mu(\a(b))=b_{(1)}\o b_{(2)}\o \a(b_{(3)})-b_{(1)}\o g\{g^{-1}, b_{(2)}\}\o b_{(3)}+g^{-1}\{g, b_{(1)}\}\o b_{(2)}\o b_{(3)},\\
&&\mu(\delta(b))=\delta(b_{(1)})\o b_{(2)}\o b_{(3)}+gb_{(1)}\o \{g^{-1}x,b_{(2)}\}\o b_{(3)}+gb_{(1)}\o g^{-1}b_{(2)}\o \delta(b_{(3)}).
\end{eqnarray*}

By the above second and third equations, we can get
$$b_{(1)}\o b_{(2)}\o \a(b_{(3)})=b_{(1)}\o g\{g^{-1},b_{(2)}\}\o b_{(3)},$$
so, by Eq. (2.1) and the above second equation, we have
$$\a(b)=g^{-1}\{g,b\}.$$

Hence, we obtain that
\begin{eqnarray*}
\a(\{a,b\})&=&g^{-1}\{g,\{a,b\}\}=-g^{-1}\{a,\{b,g\}\}-g^{-1}\{b,\{g,a\}\},\\
\{\a(a),b\}+\{a,\a(b)\}&=&\{g^{-1}\{g,a\},b\}+\{a,g^{-1}\{g,b\}\}\\
&=&-\{b,g^{-1}\{g,a\}\}+\{a,g^{-1}\{g,b\}\}\\
&=&-\{b,g^{-1}\}\{g,a\}-g^{-1}\{b,\{g,a\}\}+\{a,g^{-1}\}\{g,b\}-g^{-1}\{a,\{b,g\}\}.
\end{eqnarray*}

 Again by Eq. (4.2), we get
 $$
 \{b,g^{-1}\}\{g,a\}=\{a,g^{-1}\}\{g,b\},
 $$
 that is, $\{g^{-1}, b\}\a(a)=\{g^{-1}, a\}\a(b)$.

Conversely, if conditions (1) and (2) hold, it is easy to see $\mu(\{x,b\})=\{\mu(x),\mu(b)\})$ by a direct computation. Hence $B[x;\a,\delta]_{p}$ has a Poisson Hopf-Galois algebra structure with $B$ as a Poisson Hopf-Galois subalgebra. \ \ \ $\Box$

\vspace{4mm}

\begin{center}
{\bf \S5\quad Hopf-Galois structure on Poisson enveloping algebras}
\end{center}

\vspace{4mm}

In section, we mainly give a necessary and sufficient condition for the Poisson enveloping algebra of a Poisson Hopf-Galois algebra to be a Hopf-Galois algebra.

\vspace{2mm}

\textbf{Definition 5.1} Let $A$ be a Poisson algebra, and $(U(A),\a,\b)$ a triple, where $U(A)$ is an algebra, $\a: A\rightarrow U(A)$ is an algebra homomorphism, and $\b: A\rightarrow U(A)_{L}$ is a Lie homomorphism such that
$$
\a(\{a,b\})=[\b(a),\a(b)],~~\b(ab)=\a(a)\b(b)+\a(b)\b(a), \eqno(5.1)
$$
for all $a,b\in A$, is called the Poisson enveloping algebra in \cite{Sei2} for $A$ if $(U(A),\a,\b)$ satisfies the following: for any given algebra $B$, if $\gamma$ is an algebra homomorphism from $A$ into $B$, and $\delta$ is a Lie homomorphism from $(A,\{\cdot,\cdot\})$ into $B_{L}$ such that
$$
\gamma(\{a,b\})=[\delta(a),\gamma(b)],~~\delta(ab)=\gamma(a)\delta(b)+\gamma(b)\delta(a), \eqno(5.2)
$$
for all $a,b\in A$, then, there exists a unique algebra homomorphism $h$ from $U(A)$ into $B$ such that $h\a=\g$ and $h\b=\delta$, that is, we have the following commutative diagrams:

\vspace{4mm}

\unitlength 1mm 
\linethickness{0.4pt}
\ifx\plotpoint\undefined\newsavebox{\plotpoint}\fi 
\begin{picture}(30,28.25)(0,0)
\put(50,15.5){$A$}
\put(82,15.5){$A$}
\put(64.5,26.25){$U(A)$}
\put(66.5,3){$B$}
\put(82,18.5){\vector(-1,1){8.25}}
\put(82,14.5){\vector(-1,-1){9.75}}
\put(67.5,24){\vector(0,-1){16.25}}
\put(53,16.5){\vector(1,-1){12.25}}
\put(52.75,17.25){\vector(1,1){9.75}}
\put(54.75,22.75){$\alpha$}
\put(80.5,22.75){$\beta$}
\put(53.5,9.25){$\gamma$}
\put(78,6.25){$\delta$}
\put(68,16.25){$h$}
\end{picture}

Here $B_L$ denotes the induced Lie algebra by $B$.

\vspace{2mm}

\textbf{Remark 5.2} (1) Note that there exists a unique Poisson enveloping algebra $(U(A),\a,\b)$ up to isomorphic, that $U(A)$ is generated by $\a(A)$ and $\b(A)$ (see \cite{Sei2}, proof of Theorem 5) and that $\b(1)=0$ by Eq. (5.1), for every Poisson algebra $A$.

(2) Let $(U(A),\a,\b)$ be the Poisson enveloping algebra for a Poisson algebra $A$. In what follows, we call the two maps $\a,\b$ Poisson enveloping maps.

\vspace{2mm}

\textbf{Lemma 5.3} Suppose that $\g$ and $\d$ are two linear maps from a Poisson algebra $A$ into an algebra $B$ such that
$$
\gamma(\{a,b\})=[\delta(a),\gamma(b)],~~\delta(ab)=\gamma(a)\delta(b)+\gamma(b)\delta(a),
$$
for all $a,b\in A$, then, by \cite{Sei1},
$$
\gamma(\{a,b\})=[\g(a),\d(b)],~~\delta(ab)=\d(a)\g(b)+\d(b)\g(a). \eqno(5.3)
$$

Here the Lie structure of algebra $B$ is given by $[x,y]=xy-yx$ for $x,y\in B$.

\vspace{2mm}

\textbf{Remark 5.4} (1) By Lemma 5.3, for every Poisson enveloping algebra $(U(A),\a,\b)$ of a Poisson algebra $A$, Eq. (5.1) can also be written by
$$
\a(\{a,b\})=[\a(a),\b(b)],~~\b(ab)=\b(a)\a(b)+\b(b)\a(a).\eqno(5.1')
$$

(2) Let $A=(A,\cdot,\{\cdot,\cdot\})$ be a Poisson algebra. Define a bilinear map $\{\cdot,\cdot\}_1$ on $A$ by
$$
\{a,b\}_1=\{b,a\}
$$
for all $a,b\in A$. Then, by \cite{Sei1}, $A_1=(A,\cdot,\{\cdot,\cdot\}_1)$ is a Poisson algebra. Thus, according to Proposition 9 in \cite{Sei1}, $(U(A)^{op},\a,\b)$ is also the Poisson enveloping algebra for $A_1$, where $(U(A)^{op},\circ)$ denotes the opposite algebra of $U(A)$.

(3) Let $(U(A),\a,\b)$ be the Poisson enveloping algebra for a Poisson algebra $A$. Then, by (2), we know $\b: A\rightarrow U(A)_{L}^{op}$ is a Lie anti-homomorphism, that is, $\b(\{a,b\})=\b(b)\circ \b(a)- \b(a)\circ \b(b)$, and $\a: A\rightarrow U(A)^{op}$ is an algebra homomorphism.

Again by (5.1) and $(5.1^\prime)$, we have
$$\b(ab)=\a(a)\circ \b(b)+\a(b)\circ\b(a)=\b(a)\circ \a(b)+\b(b)\circ \a(a).\eqno(5.4)
$$
$$\a(\{a,b\})=\a(b)\circ\b(a)-\b(a)\circ\a(b)=\b(b)\circ\a(a)-\a(a)\circ\b(b).\eqno(5.5)
$$

\textbf{Lemma 5.5} Let $(U(A),\a,\b)$ be the Poisson enveloping algebra for a Poisson algebra $A$. Then

(1) $\a\o \a \o \a: A\o A\o A\longrightarrow U(A)\o U(A)^{op}\o U(A)$ is an algebra homomorphism.

Here $A\o A\o A$ and $U(A)\o U(A)^{op}\o U(A)$ denote the usual tensor product algebras.

(2) $\xi: A\o A\o A\longrightarrow (U(A)\o U(A)^{op}\o U(A))_L$ is a Lie homomorphism.

Here $\xi=(id\otimes \tau+id\otimes id\otimes id+\tau\otimes id)(\a\o \b\o \a)$, that is, $\xi=\a\o \a\o \b+\a\o \b\o \a+\b\o \a\o \a$, and the Lie structure of $A\o A\o A$ is given by Eq. (3.5).

{\bf Proof.} (1) Since $(U(A),\a,\b)$ is the Poisson enveloping algebra for a Poisson algebra $A$, the map $\a$ from $A$ to $U(A)$ is an algebra map. Hence $\a\o \a \o \a: A\o A\o A\longrightarrow U(A)\o U(A)^{op}\o U(A)$ is an algebra homomorphism.

(2) In order to prove $\xi$ to be a Lie homomorphism, we need to prove the following equality
$$
\xi(\{a\o b\o c, a'\o b'\o c'\})=[\xi(a\o b\o c),\xi(a'\o b'\o c')] \eqno(5.6)
$$
for any $a,b,c,a',b',c'\in A$.

As a matter of fact, by Remark 5.4, we get
\begin{eqnarray*}
&&\xi(\{a\o b\o c, a'\o b'\o c'\})\\
&=&\xi(\{a,a'\}\o bb'\o cc'-aa'\o \{b,b'\}\o cc'+aa'\o bb'\o \{c,c'\})\\
&=&(\a\o \a\o \b)(\{a,a'\}\o bb'\o cc'-aa'\o \{b,b'\}\o cc'+aa'\o bb'\o \{c,c'\})\\
&&+(\a\o \b\o \a)(\{a,a'\}\o bb'\o cc'-aa'\o \{b,b'\}\o cc'+aa'\o bb'\o \{c,c'\})\\
&&+(\b\o \a\o \a)(\{a,a'\}\o bb'\o cc'-aa'\o \{b,b'\}\o cc'+aa'\o bb'\o \{c,c'\})\\
&=&\a(\{a,a'\})\o \a(bb')\o \b(cc')-\a(aa')\o \a(\{b,b'\})\o \b(cc')\\
&&+\a(aa')\o \a(bb')\o \b(\{c,c'\})+\a(\{a,a'\})\o \b(bb')\o \a(cc')\\
&&-\a(aa')\o \b(\{b,b'\})\o \a(cc')+\a(aa')\o \b(bb')\o \a(\{c,c'\})\\
&&+\b(\{a,a'\})\o \a(bb')\o \a(cc')-\b(aa')\o \a(\{b,b'\})\o \a(cc')\\
&&+\b(aa')\o \a(bb')\o \a(\{c,c'\})\\
&=&\a(\{a,a'\})\o \a(bb')\o \b(cc')-\a(aa')\o \a(\{b,b'\})\o \b(cc')\\
&&+\a(aa')\o \a(bb')\o (\b(c')\circ\b(c)-\b(c)\circ\b(c'))+\a(\{a,a'\})\o \b(bb')\o \a(cc')\\
&&-\a(aa')\o (\b(b')\circ\b(b)-\b(b)\circ\b(b'))\o \a(cc')+\a(aa')\o \b(bb')\o \a(\{c,c'\})\\
&&+(\b(a')\circ\b(a)-\b(a)\circ\b(a'))\o \a(bb')\o \a(cc')-\b(aa')\o \a(\{b,b'\})\o \a(cc')\\
&&+\b(aa')\o \a(bb')\o \a(\{c,c'\})\\
&=&\underbrace{\a(\{a,a'\})\o \a(bb')\o \b(cc')}-\underbrace{\a(aa')\o \a(\{b,b'\})\o \b(cc')}\\
&&+\a(aa')\o \a(bb')\o \b(c')\circ\b(c)-\a(aa')\o \a(bb')\o\b(c)\circ\b(c')\\
&&+\underbrace{\a(\{a,a'\})\o \b(bb')\o \a(cc')}-\a(aa')\o \b(b')\circ\b(b)\o \a(cc')\\
&&+\a(aa')\o \b(b)\circ\b(b')\o \a(cc')+\underbrace{\a(aa')\o \b(bb')\o \a(\{c,c'\})}\\
&&+\b(a')\circ\b(a)\o \a(bb')\o \a(cc')-\b(a)\circ\b(a')\o \a(bb')\o \a(cc')\\
&&-\underbrace{\b(aa')\o \a(\{b,b'\})\o \a(cc')}+\underbrace{\b(aa')\o \a(bb')\o \a(\{c,c'\})}\\
&\stackrel{(5.1)(5.1')}{=}&\underbrace{\a(\{a,a'\})\o \a(bb')\o \b(c)\a(c')+\a(\{a,a'\})\o \a(bb')\o \b(c')\a(c)}\\
&&\underbrace{-\a(aa')\o\a(\{b,b'\})\o \a(c)\b(c')-\a(aa')\o\a(\{b,b'\})\o \a(c')\b(c)}
\end{eqnarray*}
\begin{eqnarray*}
&&+\a(aa')\o \a(bb')\o \b(c')\circ\b(c)-\a(aa')\o \a(bb')\o\b(c)\circ\b(c')\\
&&\underbrace{+\a(\{a,a'\})\o \b(b)\a(b')\o \a(cc')+\a(\{a,a'\})\o \b(b')\a(b)\o \a(cc')}\\
&&-\a(aa')\o \b(b')\circ\b(b)\o \a(cc')+\a(aa')\o \b(b)\circ\b(b')\o \a(cc')\\
&&\underbrace{+\a(aa')\o \b(b)\a(b')\o \a(\{c,c'\})+\a(aa')\o \b(b')\a(b)\o \a(\{c,c'\})}\\
&&+\b(a')\circ\b(a)\o \a(bb')\o \a(cc')-\b(a)\circ\b(a')\o \a(bb')\o \a(cc')\\
&&\underbrace{-\a(a)\b(a')\o \a(\{b,b'\})\o \a(cc')-\a(a')\b(a)\o \a(\{b,b'\})\o \a(cc')}\\
&&\underbrace{+\b(a)\a(a')\o \a(bb')\o \a(\{c,c'\})+\b(a')\a(a)\o \a(bb')\o \a(\{c,c'\})}\\
&\stackrel{(5.1)(5.1')}{=}&\underbrace{(\a(a)\b(a')-\b(a')\a(a))}\o \a(bb')\o \b(c)\a(c')\\
&&+\underbrace{(\b(a)\a(a')-\a(a')\b(a))}\o \a(bb')\o \b(c')\a(c)\\
&&-\a(aa')\o\underbrace{(\b(b)\a(b')-\a(b')\b(b))}\o \a(c)\b(c')\\
&&-\a(aa')\o\underbrace{(\a(b)\b(b')-\b(b')\a(b))}\o \a(c')\b(c)\\
&&+\a(aa')\o \a(bb')\o \b(c')\circ\b(c)-\a(aa')\o \a(bb')\o\b(c)\circ\b(c')\\
&&+\underbrace{(\a(a)\b(a')-\b(a')\a(a))}\o \b(b)\a(b')\o \a(cc')\\
&&+\underbrace{(\b(a)\a(a')-\a(a')\b(a))}\o \b(b')\a(b)\o \a(cc')\\
&&-\a(aa')\o \b(b')\circ\b(b)\o \a(cc')+\a(aa')\o \b(b)\circ\b(b')\o \a(cc')\\
&&+\a(aa')\o \b(b)\a(b')\o \underbrace{(\a(c)\b(c')-\b(c')\a(c))}\\
&&+\a(aa')\o \b(b')\a(b)\o \underbrace{(\b(c)\a(c')-\a(c')\b(c))}\\
&&+\b(a')\circ\b(a)\o \a(bb')\o \a(cc')-\b(a)\circ\b(a')\o \a(bb')\o \a(cc')\\
&&-\a(a)\b(a')\o \underbrace{(\b(b)\a(b')-\a(b')\b(b))}\o \a(cc')\\
&&-\a(a')\b(a)\o \underbrace{(\a(b)\b(b')-\b(b')\a(b))}\o \a(cc')\\
&&+\b(a)\a(a')\o \a(bb')\o \underbrace{(\a(c)\b(c')-\b(c')\a(c))}\\
&&+\b(a')\a(a)\o \a(bb')\o \underbrace{(\b(c)\a(c')-\a(c')}\b(c))\\
&=&\a(a)\b(a')\o \a(bb')\o \b(c)\a(c')-\a(a')\b(a)\o \a(bb')\o \b(c')\a(c)\\
&&+\a(aa')\o\b(b)\circ\a(b')\o \a(c)\b(c')-\a(aa')\o\b(b')\circ\a(b)\o \a(c')\b(c)\\
&&+\a(aa')\o \a(bb')\o \b(c)\b(c')-\a(aa')\o \a(bb')\o\b(c')\b(c)\\
&&-\b(a')\a(a)\o \a(b')\circ\b(b)\o \a(cc')+\b(a)\a(a')\o \a(b)\circ\b(b')\o \a(cc')\\
&&-\a(aa')\o \b(b')\circ\b(b)\o \a(cc')+\a(aa')\o \b(b)\circ\b(b')\o \a(cc')\\
&&-\a(aa')\o \a(b')\circ\b(b)\o \b(c')\a(c)+\a(aa')\o \a(b)\circ\b(b')\o \b(c)\a(c')\\
&&+\b(a)\b(a')\o \a(bb')\o \a(cc')-\b(a')\b(a)\o \a(bb')\o \a(cc')\\
&&+\a(a)\b(a')\o \b(b)\circ\a(b')\o \a(cc')-\a(a')\b(a)\o \b(b')\circ\a(b)\o \a(cc')\\
&&+\b(a)\a(a')\o \a(bb')\o \a(c)\b(c')-\b(a')\a(a)\o \a(bb')\o \a(c')\b(c),\\
&&[\xi(a\o b\o c),\xi(a'\o b'\o c')]
\end{eqnarray*}
\begin{eqnarray*}
&=&\xi(a\o b\o c)\xi(a'\o b'\o c')-\xi(a'\o b'\o c')\xi(a\o b\o c)\\
&=&(\a(a)\o \a(b)\o \b(c)+\a(a)\o \b(b)\o \a(c)+\b(a)\o \a(b)\o \a(c))\\
&&(\a(a')\o \a(b')\o \b(c')+\a(a')\o \b(b')\o \a(c')+\b(a')\o \a(b')\o \a(c'))\\
&&-(\a(a')\o \a(b')\o \b(c')+\a(a')\o \b(b')\o \a(c')+\b(a')\o \a(b')\o \a(c'))\\
&&(\a(a)\o \a(b)\o \b(c)+\a(a)\o \b(b)\o \a(c)+\b(a)\o \a(b)\o \a(c))\\
&=&\a(aa')\o \a(bb')\o \b(c)\b(c')+\a(aa')\o \a(b)\circ \b(b')\o \b(c)\a(c')\\
&&+\a(a)\b(a')\o \a(bb')\o \b(c)\a(c')+\a(aa')\o \b(b)\circ\a(b')\o \a(c)\b(c')\\
&&+\a(aa')\o \b(b)\circ \b(b')\o \a(cc')+\a(a)\b(a')\o \b(b)\circ\a(b')\o \a(cc')\\
&&+\b(a)\a(a')\o \a(bb')\o \a(c)\b(c')+\b(a)\a(a')\o \a(b)\circ\b(b')\o \a(cc')\\
&&+\b(a)\b(a')\o \a(bb')\o \a(cc')-\a(a'a)\o \a(b'b)\o \b(c')\b(c)\\
&&-\a(a'a)\o \a(b')\circ \b(b)\o \b(c')\a(c)-\a(a')\b(a)\o \a(b'b)\o \b(c')\a(c)\\
&&-\a(a'a)\o \b(b')\circ \a(b)\o \a(c')\b(c)-\a(a'a)\o \b(b')\circ \b(b)\o \a(c'c)\\
&&-\a(a')\b(a)\o \b(b')\circ \a(b)\o \a(c'c)-\b(a')\a(a)\o \a(b'b)\o \a(c')\b(c)\\
&&-\b(a')\a(a)\o \a(b')\circ \b(b)\o \a(c'c)-\b(a')\b(a)\o \a(b'b)\o \a(c'c).
\end{eqnarray*}

By the above proof, we know that $\xi(\{a\o b\o c, a'\o b'\o c'\})=[\xi(a\o b\o c),\xi(a'\o b'\o c')]$. So, $\xi$ is a Lie homomorphism.\ \ \ $\Box$

\vspace{2mm}

\textbf{Lemma 5.6} Let $A$ and $B$ be Poisson algebras and $C$ an algebra. If $\phi: A\rightarrow B$ is a Poisson homomorphism, $\a: B\rightarrow C$ an algebra homomorphism and $\b: B\rightarrow C_L$ a Lie homomorphism such that
$$
\a(\{b,b'\})=[\b(b),\a(b')],~~\b(bb')=\a(b)\b(b')+\a(b')\b(b)
$$
for all $b,b'\in B$, then, according to \cite{Sei1}, $\a\phi: A\rightarrow C$ is an algebra homomorphism and $\b\phi: A\rightarrow C_L$ a Lie homomorphism such that
$$
\a\phi(\{a,a'\})=[\b\phi(a),\a\phi(a')],~~\b\phi(aa')=\a\phi(a)\b\phi(a')+\a\phi(a')\b\phi(a)
$$
for all $a,a'\in A$.

\vspace{2mm}

\textbf{Lemma 5.7} Let $(U(A),\a, \b)$ be the Poisson enveloping algebra for a Poisson algebra $A$. Then $(U(A)\o U(A)^{op}\o U(A), \a\o \a\o \a, \xi)$ is the Poisson enveloping algebra for $A\o A\o A$.

 Here $\xi=(id\otimes \tau+id\otimes id\otimes id+\tau\otimes id)(\a\o \b\o \a)$, that is, $\xi=\a\o \a\o \b+\a\o \b\o \a+\b\o \a\o \a$.

{\bf Proof.} In order to prove that $(U(A)\o U(A)^{op}\o U(A), \a\o \a\o \a, \xi)$ is the Poisson enveloping algebra for $A\o A\o A$, the Proof is divided into the following steps by Definition 5.1.

 $\bullet$ Step 1, we prove that $\a\o \a\o \a$ is an algebra homomorphism and $\xi$ is a Lie homomorphism.

 It is obvious by Lemma 5.5.

 $\bullet$ Step 2, we prove the equalities $(5.1)$ hold for the Poisson enveloping map $\a\o \a\o \a$.

For any $a,b,c,a',b',c'\in A$, we get
\begin{eqnarray*}
&&(\a\o \a \o \a)(\{a\o b\o c, a'\o b'\o c'\})\\
&=&(\a\o \a \o \a)(\{a,a'\}\o bb'\o cc'-aa'\o \{b,b'\}\o cc'+aa'\o bb'\o \{c,c'\})\\
\end{eqnarray*}
\begin{eqnarray*}
&=&\a(\{a,a'\})\o \a(bb')\o \a(cc')-\a(aa')\o \a(\{b,b'\})\o \a(cc')+\a(aa')\o \a(bb')\o \a(\{c,c'\})\\
&\stackrel{(5.1)}{=}&(\underbrace{\b(a)\a(a')-\a(a')\b(a)})\o \a(bb')\o \a(cc')-\a(aa')\o (\underbrace{\a(b')\circ\b(b)-\b(b)\circ\a(b')})\o \a(cc')\\
&&+\a(aa')\o \a(bb')\o (\underbrace{\b(c)\a(c')-\a(c')\b(c)})\\
&=&\b(a)\a(a')\o \a(bb')\o \a(cc')-\a(a')\b(a)\o \a(bb')\o \a(cc')-\a(aa')\o\a(b')\circ\b(b)\o \a(cc')\\
&&+\a(aa')\o \b(b)\circ\a(b')\o \a(cc')+\a(aa')\o \a(bb')\o \b(c)\a(c')-\a(aa')\o \a(bb')\o \a(c')\b(c),\\
&&[\xi(a\o b\o c), (\a\o \a\o \a)(a'\o b'\o c')]\\
&=&\xi(a\o b\o c)(\a(a')\o \a(b')\o \a(c'))-(\a(a')\o \a(b')\o \a(c'))\xi(a\o b\o c)\\
&=&(\underbrace{\a(a)\o \a(b)\o \b(c)+\a(a)\o \b(b)\o \a(c)+\b(a)\o \a(b)\o \a(c)})(\a(a')\o \a(b')\o \a(c'))\\
&&-(\a(a')\o \a(b')\o \a(c'))(\underbrace{\a(a)\o \a(b)\o \b(c)+\a(a)\o \b(b)\o \a(c)+\b(a)\o \a(b)\o \a(c)})\\
&=&\a(aa')\o \a(bb')\o \b(c)\a(c')+\a(aa')\o \b(b)\circ \a(b')\o \a(cc')\\
&&+\b(a)\a(a')\o \a(bb')\o \a(cc')-\a(aa')\o \a(bb')\o \a(c')\b(c)\\
&&-\a(aa')\o \a(b')\circ \b(b)\o \a(cc')-\a(a')\b(a)\o \a(bb')\o \a(cc'),
\end{eqnarray*}
so, $(\a\o \a \o \a)(\{a\o b\o c, a'\o b'\o c'\})=[\xi(a\o b\o c), (\a\o \a\o \a)(a'\o b'\o c')]$.
\begin{eqnarray*}
&&\xi(aa'\o bb'\o cc')\\
&=&\a(aa')\o \a(bb')\o \b(cc')+\a(aa')\o \b(bb')\o \a(cc')+\b(aa')\o \a(bb')\o \a(cc')\\
&\stackrel{(5.1)(5.1')}{=}&\a(aa')\o \a(bb')\o (\underbrace{\a(c)\b(c')+\a(c')\b(c)})+\a(aa')\o (\underbrace{\a(b')\circ\b(b)+\a(b)\circ\b(b')})\o \a(cc')\\
&&+(\underbrace{\a(a)\b(a')+\a(a')\b(a)})\o \a(bb')\o \a(cc'),\\
&&(\a\o \a \o \a)(a\o b\o c)\xi(a'\o b'\o c')+(\a\o \a \o \a)(a'\o b'\o c')\xi(a\o b\o c)\\
&=&(\a(a)\o \a(b)\o \a(c))(\underbrace{\a(a')\o \a(b')\o \b(c')+\a(a')\o \b(b')\o \a(c')+\b(a')\o \a(b')\o \a(c')})\\
&&+(\a(a')\o \a(b')\o \a(c'))(\underbrace{\a(a)\o \a(b)\o \b(c)+\a(a)\o \b(b)\o \a(c)+\b(a)\o \a(b)\o \a(c)})\\
&=&\a(aa')\o \a(bb')\o \a(c)\b(c')+\a(aa')\o \a(b)\circ\b(b')\o \a(cc')+\a(a)\b(a')\o \a(bb')\o \a(cc')\\
&&+\a(aa')\o \a(bb')\o \a(c')\b(c)+\a(aa')\o \a(b')\circ\b(b)\o \a(cc')+\a(a')\b(a)\o \a(bb')\o \a(cc'),
\end{eqnarray*}
so, $\xi(aa'\o bb'\o cc')=(\a\o \a \o \a)(a\o b\o c)\xi(a'\o b'\o c')+(\a\o \a \o \a)(a'\o b'\o c')\xi(a\o b\o c)$.

 $\bullet$ Step 3, we prove the universal property of $(U(A)\o U(A)^{op}\o U(A), \a\o \a\o \a, \xi)$.

Given an algebra $B$, let $m_B$ be the multiplication map on $B$. If $\g$ is an algebra homomorphism from $A\o A\o A$ into $B$ and $\d$ is a Lie homomorphism from $A\o A\o A$ into $B_L$ such that
 $$
\g(\{a\o b\o c, a'\o b'\o c'\})=[\d(a\o b\o c),\g(a'\o b'\o c')],\eqno{(5.7)}
$$
$$
\d((a\o b\o c)(a'\o b'\o c'))=\g(a\o b\o c)\d(a'\o b'\o c')+\g(a'\o b'\o c')\d(a\o b\o c)\eqno{(5.8)}
$$
 for all $a,b,c,a',b',c'\in A$.

 Let $i_1, i_2$ and $i_3$ be the Poisson homomorphisms defined by
 \begin{eqnarray*}
 &&i_1:A\rightarrow A\o A\o A,~~i_1(a)=a\o 1\o 1,\\
 &&i_2:A_1\rightarrow A\o A\o A,~i_2(a)=1\o a\o 1,\\
 &&i_3:A\rightarrow A\o A\o A,~~i_3(a)=1\o 1\o a
 \end{eqnarray*}
 for all $a\in A$.

 Then, by Lemma 5.6 and the universal property of a Poisson algebra $A$, there exist algebra homomorphisms $f,\theta:U(A)\rightarrow B$, and $g:U(A)^{op}\rightarrow B$ such that $f\a=\g i_1,~f\b=\d i_1,~g\a=\g i_2,~g\b=\d i_2,~\theta\a=\g i_3,~\theta\b=\d i_3$, that is, we have the following commutative diagrams:

\unitlength 1mm 
\linethickness{0.4pt}
\ifx\plotpoint\undefined\newsavebox{\plotpoint}\fi 
\begin{picture}(142.25,31.25)(0,0)
\put(2,3.5){$A$}
\put(1.5,20){$U(A)$}
\put(34,20){$B$}
\put(9,4.25){\vector(1,0){16}}
\put(12,21){\vector(1,0){20}}
\put(4.25,7.25){\vector(0,1){11.25}}
\put(35,7.25){\vector(0,1){11.25}}
\put(25.75,3.5){$A\otimes A\otimes A$}
\put(16,5.5){$i_1$}
\put(4.5,12){$\alpha,\beta$}
\put(18,22){$f$}
\put(29.5,12){$\gamma,\delta$}
\put(52,3.75){$A_1$}
\put(76.25,3.75){$A\otimes A\otimes A$}
\put(50,20){$U(A)^{op}$}
\put(84,20){$B$}
\put(64,21.5){\vector(1,0){20}}
\put(53.75,7.25){\vector(0,1){11.25}}
\put(84.5,7.25){\vector(0,1){11.25}}
\put(58,4.75){\vector(1,0){18.5}}
\put(53.75,12.75){$\alpha,\beta$}
\put(78.5,12.75){$\gamma,\delta$}
\put(64.5,6.25){$i_2$}
\put(73,22.25){$g$}
\put(106.5,3.5){$A$}
\put(105.25,20){$U(A)$}
\put(108,7.25){\vector(0,1){11.25}}
\put(135.5,7.25){\vector(0,1){11.25}}
\put(111.25,4.5){\vector(1,0){14.25}}
\put(114.25,21.5){\vector(1,0){18}}
\put(127,3.5){$A\otimes A\otimes A$}
\put(134,20){$B$}
\put(108,12.25){$\alpha,\beta$}
\put(130.25,12.25){$\gamma,\delta$}
\put(116.5,6.75){$i_3$}
\put(121.75,22.25){$\theta$}
\end{picture}

 Moreover, we have
$$\d i_1(a)\g i_2(b)=\g i_2(b)\d i_1(a),\eqno{(5.9)}
$$
$$\d(a\o b\o 1)\g(1\o 1\o c)=\g(1\o 1\o c)\d(a\o b\o 1)\eqno{(5.10)}
$$
for all $a,b\in A$. This is because
 \begin{center}
$[\d i_1(a), \g i_2(b)]\stackrel{(5.7)}=\g(\{a\o 1\o 1,1\o b\o 1\})\stackrel{(3.5)}=0,$
 \end{center}
  \begin{center}
 $[\d(a\o b\o 1), \g (1\o 1\o c)]\stackrel{(5.7)}=\g(\{a\o b\o 1,1\o 1\o c\})\stackrel{(3.5)}=0.$
 \end{center}

 In a similar way, by Eq. $(5.7)$, we can prove that
 $$\d i_m(a)\g i_n(b)=\g i_n(b)\d i_m(a)\eqno{(5.9^\prime)}
 $$
 for any $1\leqslant m\neq n\leqslant 3$.

 Because $\d$ is a Lie homomorphism from $A\o A\o A$ to $B$, we can directly prove that
 $$\d i_m(a)\d i_n(b)=\d i_n(b)\d i_m(a)\eqno{(5.11)}
 $$
 for any $1\leqslant m\neq n\leqslant 3$.

 Hence, for all $a,b,c\in A$, we have
\begin{eqnarray*}
&& m_B(id\o m_B)(f\o g\o \theta)(\a\o \a\o \a)(a\o b\o c)\\
 &=&f(\a(a))g(\a(b))\theta(\a(c))=\g i_1(a) \g i_2(a) \g i_3(a)\\
 &=&\g(i_1(a)i_2(b)i_3(c))=\g(a\o b\o c),\\
 &&m_B(id\o m_B)(f\o g\o \theta)\xi(a\o b\o c)\\
&=&f(\a(a))g(\a(b))\theta(\b(c))+f(\a(a))g(\b(b))\theta(\a(c))+f(\b(a))g(\a(b))\theta(\a(c))\\
&=&\g i_1(a)\g i_2(b)\d i_3(c)+\g i_1(a)\d i_2(b)\g i_3(c)+\d i_1(a)\g i_2(b)\g i_3(c)\\
&\stackrel{(5.9)}=&\g i_1(a)\g i_2(b)\d i_3(c)+\g i_1(a)\d i_2(b)\g i_3(c)+\g i_2(b)\d i_1(a)\g i_3(c)\\
&=&\g i_1(a)\g i_2(b)\d i_3(c)+\underbrace{(\g i_1(a)\d i_2(b)+\g i_2(b)\d i_1(a))}\g i_3(c)\\
&\stackrel{(5.8)}{=}&\g i_1(a)\g i_2(b)\d i_3(c)+\d(i_1(a)i_2(b))\g i_3(c)\\
&=&\g(a\o b\o 1)\d (1\o 1\o c)+\d (a\o b\o 1)\g (1\o 1\o c)\\
&\stackrel{(5.10)}=&\g(a\o b\o 1)\d (1\o 1\o c)+\g (1\o 1\o c)\d (a\o b\o 1)\\
&\stackrel{(5.8)}{=}&\d(a\o b\o c),
 \end{eqnarray*}
 for all $a,b,c\in A$, that is, we have
$$m_B(id\o m_B)(f\o g\o \theta)(\a\o \a\o \a)=\g, \eqno{(5.12)}
$$
$$m_B(id\o m_B)(f\o g\o \theta)\xi=\d. \eqno{(5.13)}
$$

 In what follows, we prove that $m_B(id\o m_B)(f\o g\o \theta)$ is an algebra homomorphism from $U(A)\o U(A)^{op}\o U(A)$ to $B$.

 Since $\g$ is an algebra homomorphism, we are easy to see that
 $$\g(a\o b\o c)\g(a'\o b'\o c')=\g(a'\o b'\o c')\g(a\o b\o c)\eqno{(5.14)}$$
  for any $a,b,c,a',b',c'\in A$.

 For any $a,b,c,a',b',c'\in A$, we have,
\begin{eqnarray*}
 &&m_B(id\o m_B)(f\o g\o \theta)((\a(a)\o \b(b)\o \a(c))(\a(a')\o \a(b')\o \b(c')))\\
  &=&m_B(id\o m_B)(f\o g\o \theta)(\a(a)\a(a')\o \b(b)\circ\a(b')\o \a(c)\b(c'))\\
    &=&m_B(id\o m_B)(f(\a(a)\a(a'))\o g(\b(b)\circ\a(b'))\o  \theta(\a(c)\b(c')))\\
    &=&m_B(id\o m_B)(f(\a(a))f(\a(a'))\o g(\b(b))g(\a(b'))\o  \theta(\a(c))\theta(\b(c')))\\
    &=&f(\a(a))f(\a(a'))g(\b(b))g(\a(b'))\theta(\a(c))\theta(\b(c'))\\
    &=&\g i_1(a)\underbrace{\g i_1(a')\d i_2(b)}\g i_2(b')\g i_3(c)\d i_3(c')\stackrel{(5.9^\prime)}=\g i_1(a)\d i_2(b)\g i_1(a')\underbrace{\g i_2(b')\g i_3(c)}\d i_3(c')\\
    &\stackrel{(5.14)}=&\g i_1(a)\d i_2(b)\underbrace{\g i_1(a')\g i_3(c)}\g i_2(b')\d i_3(c')=\g i_1(a)\d i_2(b)\g i_3(c)\g i_1(a')\g i_2(b')\d i_3(c')\\
    &=&m_B(id\o m_B)(f\o g\o \theta)((\a(a)\o \b(b)\o \a(c)))m_B(id\o m_B)(f\o g\o \theta)\\
    &\  &((\a(a')\o \a(b')\o \b(c'))).
    \end{eqnarray*}

 Similarly, we can prove the other cases are satisfied, so, $m_B(id\o m_B)(f\o g\o \theta)$ is an algebra homomorphism such that
 $m_B(id\o m_B)(f\o g\o \theta)(\a\o \a\o \a)=\g,~m_B(id\o m_B)(f\o g\o \theta)\xi=\d$ by the equalities (5.10) and (5.11).

 Finally, we prove that the algebra homomorphism from $U(A\o A\o A)$ to $B$ is unique.

 Suppose that there exists a map $h:U(A\o A\o A)\rightarrow B$ is an algebra homomorphism such that
 $h(\a\o \a\o \a)=\g,~h\xi=\d.$ Then, for all $a\in A$, we have
 \begin{eqnarray*}
 m_B(id\o m_B)(f\o g\o \theta)(\a(a)\o 1\o 1)&=& m_B(id\o m_B)(f\o g\o \theta)(\a\o \a\o \a)(a\o 1\o 1)\\
 &\stackrel{(5.12)}=&\g(a\o 1\o 1)=h(\a\o \a\o \a)(a\o 1\o 1)\\
 &=&h(\a(a)\o 1\o 1).
 \end{eqnarray*}

 In a similar way, we can prove that
 \begin{eqnarray*}
 m_B(id\o m_B)(f\o g\o \theta)(1\o \a(a)\o 1)&=&h(1\o \a(a)\o 1),\\
 m_B(id\o m_B)(f\o g\o \theta)(1\o 1\o \a(a))&=&h(1\o 1\o \a(a)).
  \end{eqnarray*}

 Moreover, we know that
  \begin{eqnarray*}
 m_B(id\o m_B)(f\o g\o \theta)(\b(a)\o 1\o 1)&=&m_B(id\o m_B)(f\o g\o \theta)\xi(a\o 1\o 1)\\
 &\stackrel{(5.13)}=&\d(a\o 1\o 1)=h\xi(a\o 1\o 1)\\
 &=&h(\b(a)\o 1\o 1).
 \end{eqnarray*}

 In a similar way, we can prove that
 \begin{eqnarray*}
 m_B(id\o m_B)(f\o g\o \theta)(1\o \b(a)\o 1)&=&h(1\o \b(a)\o 1),\\
 m_B(id\o m_B)(f\o g\o \theta)(1\o 1\o \b(a))&=&h(1\o 1\o \b(a)).
 \end{eqnarray*}

Hence, by the above proof and the fact that $U(A)$ is generated by $\a(A)$ and $\b(A)$, we obtain that $m_B(id\o m_B)(f\o g\o \theta)=h$, which completes the proof. \ \ \  $\Box$

\vspace{2mm}

The following lemma can be found in \cite{Sei1}.

\vspace{2mm}

\textbf{Lemma 5.8} Let $(U(A),\a_A,\b_A)$ and $(U(B),\a_B,\b_B)$ be Poisson enveloping algebras for Poisson algebras $A$ and $B$, respectively. If $\phi:A\rightarrow B$ is a Poisson homomorphism, then, there exists a unique algebra homomorphism $U(\phi):U(A)\rightarrow U(B)$ such that
$$U(\phi)\a_A=\a_B\phi, \ \  U(\phi)\b_A=\b_B\phi,\eqno{(5.15)}
$$
that is, the following diagram is commutative:

\vspace{3mm}

\unitlength 1mm 
\linethickness{0.4pt}
\ifx\plotpoint\undefined\newsavebox{\plotpoint}\fi 
\begin{picture}(50.75,24.75)(0,0)
\put(52,22.25){$U(A)$}
\put(53,2.75){$A$}
\put(77,22.5){$U(B)$}
\put(78,2.75){$B$}
\put(58,3.85){\vector(1,0){18.5}}
\put(60.5,23.75){\vector(1,0){15.5}}
\put(54.75,7){\vector(0,1){14}}
\put(79.25,6.5){\vector(0,1){14}}
\put(43.5,13.5){$\alpha_A,\beta_A$}
\put(80.75,13.5){$\alpha_B,\beta_B$}
\put(63.25,24.75){$U(\phi)$}
\put(65.25,5.5){$\phi$}
\end{picture}

\vspace{2mm}

According to the above lemmas, we obtain a main result in this section.

\vspace{2mm}

\textbf{Theorem 5.9} Let $(A,\mu)$ be a Poisson Hopf-Galois algebra, and $(U(A),\a,\b)$ its Poisson enveloping algebra. Then
$(U(A),U(\mu))$ satisfies the following conditions:
$$
U(\mu)\a=(\a\o \a\o \a)\mu,\ \ U(\mu)\b=\xi\mu.\eqno{(5.16)}
$$

 Here $\xi=(id\otimes \tau+id\otimes id\otimes id+\tau\otimes id)(\a\o \b\o \a)$, that is, $\xi=\a\o \a\o \b+\a\o \b\o \a+\b\o \a\o \a$.

Thus, $(U(A),U(\mu))$ is a Hopf-Galois algebra if and only if
$$\b+m_{U(A)}(id\o m_{U(A)})(\a\o \b\o \a)\mu=0,\eqno{(5.17)}
$$
that is, $\b(a)=-\a(a_{(1)})\b(a_{(2)})\a(a_{(3)})$ for any $a\in A$.

{\bf Proof.} Since $\mu$ is a Poisson homomorphism and $(U(A)\o U(A)^{op}\o U(A),\a\o \a\o \a, \xi)$ is the Poisson enveloping algebra of $A\o A\o A$, there exists an algebra homomorphism
$U(\mu):U(A)\rightarrow U(A)\o U(A)^{op}\o U(A)$ by Lemma 5.8, such that the equalities (5.16) hold.

 ``$\Longleftarrow$" According to the equalities (5.16), we have
$$ \alpha(a)_{(1)}\otimes  \alpha(a)_{(2)}\otimes \alpha(a)_{(3)}=\alpha(a_{(1)})\otimes  \alpha(a_{(2)})\otimes \alpha(a_{(3)}),\eqno{(5.18)}
$$
$$
  \beta(a)_{(1)}\otimes  \beta(a)_{(2)}\otimes \beta(a)_{(3)}=\xi(a_{(1)})\otimes  \xi(a_{(2)})\otimes \xi(a_{(3)})\eqno{(5.19)}
  $$
for any $a\in A$, where $U(\mu)(\alpha(a))$ is denoted by $\alpha(a)_{(1)}\otimes  \alpha(a)_{(2)}\otimes \alpha(a)_{(3)}$, and $U(\mu)(\beta(a))$ denoted by $\beta(a)_{(1)}\otimes  \beta(a)_{(2)}\otimes \beta(a)_{(3)}$.

By Eq. (5.18) and $\mu$ being a Hopf-Galois map for $A$, we are easy to check that
$$U(\mu)(\a(a)_{(1)})\o \a(a)_{(2)}\o\a(a)_{(3)}=\a(a)_{(1)}\o \a(a)_{(2)}\o U(\mu)(\a(a)_{(3)}).\eqno{(5.20)}
$$

In addition, for any $a\in A$, we have
\begin{eqnarray*}
&&U(\mu)(\b(a)_{(1)})\o \b(a)_{(2)}\o\b(a)_{(3)}\\
&\stackrel{(5.19)}=&U(\mu)(\a(a_{(1)}))\o\a(a_{(2)})\o\b(a_{(3)})+U(\mu)(\a(a_{(1)}))\o \b(a_{(2)})\o\a(a_{(3)})\\
&&+U(\mu)(\b(a_{(1)}))\o \a(a_{(2)})\o\a(a_{(3)})\\
&\stackrel{(5.18),(5.19)}=&\a(a_{(1)(1)})\o \a(a_{(1)(2)})\o \a(a_{(1)(3)})\o\a(a_{(2)})\o\b(a_{(3)})\\
&&+\a(a_{(1)(1)})\o \a(a_{(1)(2)})\o \a(a_{(1)(3)})\o \b(a_{(2)})\o\a(a_{(3)})\\
&&+\underbrace{\a(a_{(1)(1)})\o \a(a_{(1)(2)})\o \b(a_{(1)(3)})}\o \a(a_{(2)})\o\a(a_{(3)})\\
&&+\underbrace{\a(a_{(1)(1)})\o \b(a_{(1)(2)})\o \a(a_{(1)(3)})}\o \a(a_{(2)})\o\a(a_{(3)})\\
&&+\underbrace{\b(a_{(1)(1)})\o \a(a_{(1)(2)})\o \a(a_{(1)(3)})}\o \a(a_{(2)})\o\a(a_{(3)}),\\
&&\b(a)_{(1)}\o \b(a)_{(2)}\o U(\mu)(\b(a)_{(3)})\\
&\stackrel{(5.19)}=&\a(a_{(1)})\o\a(a_{(2)})\o U(\mu)(\b(a_{(3)}))+\a(a_{(1)})\o \b(a_{(2)})\o U(\mu)(\a(a_{(3)}))\\
&&+\b(a_{(1)})\o \a(a_{(2)})\o U(\mu)(\a(a_{(3)}))\\
&=&\a(a_{(1)})\o\a(a_{(2)})\o \underbrace{\a(a_{(3)(1)})\o \a(a_{(3)(2)})\o \b(a_{(3)(3)})}\\
&&+\a(a_{(1)})\o\a(a_{(2)})\o \underbrace{\a(a_{(3)(1)})\o \b(a_{(3)(2)})\o \a(a_{(3)(3)})}\\
&&+\a(a_{(1)})\o\a(a_{(2)})\o \underbrace{\b(a_{(3)(1)})\o \a(a_{(3)(2)})\o \a(a_{(3)(3)})}\\
&&+\a(a_{(1)})\o \b(a_{(2)})\o \underline{\a(a_{(3)(1)})\o \a(a_{(3)(2)})\o \a(a_{(3)(3)})}\\
&&+\b(a_{(1)})\o \a(a_{(2)})\o \underline{\a(a_{(3)(1)})\o \a(a_{(3)(2)})\o \a(a_{(3)(3)})}\\
&\stackrel{(2.2)}=&\a(a_{(1)(1)})\o \a(a_{(1)(2)})\o \a(a_{(1)(3)})\o\a(a_{(2)})\o\b(a_{(3)})\\
&&+\a(a_{(1)(1)})\o \a(a_{(1)(2)})\o \a(a_{(1)(3)})\o \b(a_{(2)})\o\a(a_{(3)})\\
&&+\a(a_{(1)(1)})\o \a(a_{(1)(2)})\o \b(a_{(1)(3)})\o \a(a_{(2)})\o\a(a_{(3)})\\
&&+\a(a_{(1)(1)})\o \b(a_{(1)(2)})\o \a(a_{(1)(3)})\o \a(a_{(2)})\o\a(a_{(3)})\\
&&+\b(a_{(1)(1)})\o \a(a_{(1)(2)})\o \a(a_{(1)(3)})\o \a(a_{(2)})\o\a(a_{(3)}),
\end{eqnarray*}
so, we have
$$U(\mu)(\b(a)_{(1)})\o \b(a)_{(2)}\o\b(a)_{(3)}=\b(a)_{(1)}\o \b(a)_{(2)}\o U(\mu)(\b(a)_{(3)}). \eqno{(5.21)}
$$

In what follows, we prove Eq. $(2.1)$ of an algebra map $U(\mu)$ holds.

Assume that $\b(a)=-\a(a_{(1)})\b(a_{(2)})\a(a_{(3)})$ for any $a\in A$. Then, we have
\begin{eqnarray*}
&&\b(a)_{(1)}\o \b(a)_{(2)}\b(a)_{(3)}\\
&\stackrel{(5.19)}=&\a(a_{(1)})\o \a(a_{(2)})\b(a_{(3)})+\a(a_{(1)})\o \b(a_{(2)})\a(a_{(3)})+\b(a_{(1)})\o \a(a_{(2)})\a(a_{(3)})\\
&\stackrel{(2.1)}=&\a(a_{(1)})\o\a(a_{(2)})\b(a_{(3)})+\a(a_{(1)})\o \b(a_{(2)})\a(a_{(3)})+\b(a)\o 1\\
&=&-\underbrace{\a(a_{(1)})\o \a(a_{(2)})\a(a_{(3)(1)})\b(a_{(3)(2)})\a(a_{(3)(3)})}+\a(a_{(1)})\o \b(a_{(2)})\a(a_{(3)})+\b(a)\o 1\\
&=&-\a(a_{(1)(1)})\o \a(a_{(1)(2)})\a(a_{(1)(3)})\b(a_{(2)})\a(a_{(3)})+\a(a_{(1)})\o \b(a_{(2)})\a(a_{(3)})+\b(a)\o 1
\end{eqnarray*}
\begin{eqnarray*}
&=&-\a(a_{(1)})\o \b(a_{(2)})\a(a_{(3)})+\a(a_{(1)})\o \b(a_{(2)})\a(a_{(3)})+\b(a)\o 1\\
&=&\b(a)\o 1.
\end{eqnarray*}

In a similar way, we can prove that
 \begin{eqnarray*}
 \b(a)_{(1)}\b(a)_{(2)}\o \b(a)_{(3)}=1\o \b(a)
 \end{eqnarray*}
for any $a\in A$. Hence $(U(A),U(\mu))$ is a Hopf-Galois algebra.

 ``$\Longleftarrow"$ If $(U(A),U(\mu))$ is a Hopf-Galois algebra, then $\b(a)_{(1)}\o \b(a)_{(2)}\b(a)_{(3)}=\b(a)\o 1$ for any $a\in A$.

According to Eq. (5.19), we know that
 \begin{eqnarray*}
 \b(a)_{(1)}\o \b(a)_{(2)}\o \b(a)_{(3)}&=&\a(a_{(1)})\o \a(a_{(2)})\o \b(a_{(3)})\\
 &+&\a(a_{(1)})\o \b(a_{(2)})\o \a(a_{(3)})+\b(a_{(1)})\o \a(a_{(2)})\o \a(a_{(3)}).
 \end{eqnarray*}

By applying $m_{U(A)}(id\o m_{U(A)})$ to the above equation, we get
$$\b(a)=-\a(a_{(1)})\b(a_{(2)})\a(a_{(3)}),$$
which completes this proof. \ \ \ \ $\Box$

\vspace{4mm}
\newpage

\makeatletter
\newcommand{\adjustmybblparameters}{\setlength{\itemsep}{0\baselineskip}\setlength{\parsep}{0pt}}
\let\ORIGINALlatex@openbib@code=\@openbib@code
\renewcommand{\@openbib@code}{\ORIGINALlatex@openbib@code\adjustmybblparameters}
\makeatother

\end{document}